\documentclass[twoside,11pt]{article} %
\usepackage[short,c2]{optidef}
\usepackage{tikz, lscape, amsmath,amsfonts,latexsym,amssymb, amsthm}
\usepackage[utf8]{inputenc}
\usepackage[english]{babel}
\usepackage{lmodern, verbatim}
\usepackage[T1]{fontenc}
\usepackage{calc,bm,nicefrac}
\usepackage{tabularx}
\usepackage[shortlabels]{enumitem}
\usepackage[mathscr]{eucal}
\usepackage{subfig,stackrel,floatrow,url, empheq}
\usepackage[colorinlistoftodos,bordercolor=orange,textwidth=2cm,backgroundcolor=orange!20,linecolor=orange,textsize=scriptsize]{todonotes}
\renewcommand{\b}{\mathbf}

\usepackage[round]{natbib}

\usepackage[top=3cm, bottom=3cm, left=2cm, right=2cm]{geometry}
\usepackage{hyperref}
\usepackage[font=footnotesize]{caption}
\usepackage[nameinlink,capitalize]{cleveref}
\usepackage{stackengine}
\usepackage{scalerel}
\usepackage{xcolor}
\usepackage{cite}
\usepackage{booktabs}
\usepackage{macro}

\definecolor{dgreen}{rgb}{0.00,0.49,0.00}
\definecolor{dblue}{rgb}{0,0.08,0.75}
\hypersetup{
	colorlinks = true,
	citecolor=dgreen,
	urlcolor=dblue,
	linkcolor=dblue
}
\usepackage{bibunits}
\usepackage{algorithmicx}
\usepackage{algpseudocode}

\newcommand{\V}{\mathcal{V}}
\newcommand{\U}{\mathcal{U}}

\newcommand{\BB}{\mathbb{B}}

\newcommand{\NN}{\mathbb{N}}

\newcommand{\RR}{\mathbb{R}}

\newcommand{\Hk}{\mathscr{H}_k}

\newcommand{\FK}{\mathscr{F}_K}
\newcommand{\HK}{\mathscr{H}_K}
\newcommand{\Psc}{\mathscr{P}}
\newcommand{\Lcal}{\mathcal{L}}
\newcommand{\Ecal}{\mathcal{E}}

\newcommand{\Ical}{\mathcal{I}}
\newcommand{\Ucal}{\mathcal{U}}
\newcommand{\Kcons}{\mathscr{K}}
\newcommand{\F}{\mathcal{F}}      %
\newcommand{\E}{\mathcal{E}}      %

\newcommand{\latr}{\la_{\text{tr}}}

\newcommand{\iv}[2]{[\![#1,#2]\!]} %
\renewcommand{\d}{\mathrm{d}}    %

\newcommand{\D}{\mathcal{D}}  %
\newcommand{\dBB}{\mathring{\mathbb{B}}}     %

\renewcommand{\b}{\mathbf}        %
\newcommand{\tb}{\textbf}         %
\newcommand{\p}{\partial}         %
\newcommand{\Rnn}{\R_{+}}         %
\newcommand{\Zp}{\N^*}            %

\DeclareMathOperator*{\Id}{Id}

\DeclareMathOperator*{\Sp}{span}

\DeclareMathOperator*{\Ker}{ker}
\DeclareMathOperator*{\diag}{diag} %
\DeclareMathOperator*{\Dom}{Dom}

\renewcommand{\d}{\mathrm{d}} %
\renewcommand{\b}{\mathbf}
\renewcommand{\hh}{\mathcal H}

\newcommand{\Hphi}{\mathscr{H}_\phi}

\newcommand{\epsArr}{\bm{\epsilon}}%

\newcommand{\subscript}[2]{$#1 _ #2$}

\newlist{assumplist}{enumerate}{1}
\setlist[assumplist]{label=(\subscript{\textbf{A}}{{\arabic*}})}
\Crefname{assumplisti}{Assumption}{Assumptions}

\newlist{assumplistp}{enumerate}{1}
\setlist[assumplistp]{label=(\subscript{\textbf{A}'}{{\arabic*}})}
\Crefname{assumplistpi}{Assumption}{Assumptions}

\newcommand{\rv}{\textcolor{black}}

\date{\today}

\title{Approximation of optimization problems with constraints through kernel Sum-Of-Squares}
\author{Pierre-Cyril Aubin-Frankowski, Alessandro Rudi\footnote{INRIA and Département d’Informatique, École Normale Supérieure, 
	PSL Research University, \emph{pierre-cyril.aubin@inria.fr, alessandro.rudi@inria.fr}}}

\date{\today}

\begin{document}
	\maketitle
	\begin{abstract}
	    Handling an infinite number of inequality constraints in infinite-dimensional spaces occurs in many fields, from global optimization to optimal transport. These problems have been tackled individually in several previous articles through kernel Sum-Of-Squares (kSoS) approximations. We propose here a unified theorem to prove convergence guarantees for these schemes. Pointwise inequalities are turned into equalities within a class of nonnegative kSoS functions. Assuming further that the functions appearing in the problem are smooth, focusing on pointwise equality constraints enables the use of scattering inequalities to mitigate the curse of dimensionality in sampling the constraints. Our approach is illustrated in learning vector fields with side information, here the invariance of a set.
	    
	    \tb{Keywords:} Reproducing kernels, nonconvex optimization, constraints, sum-of-squares.
	\end{abstract}
	\section{Introduction}
	
	In this study, we are interested in solving approximately the following family of convex optimization problems
	\begin{argminie}|s|
		{\b f=[f_1;\dots;f_P] \in \C^0(\X,\R^P)}{\Lcal(\b f)}{\label{opt-cons_gen}}{\bar{\b f} \in}%
		\addConstraint{\b f(x)}{\in \b F(x),\, \forall\, x\in \X \label{eq:cons_gen_F}}
		\addConstraint{\b f}{\in \E, \label{eq:cons_gen_E}}
	\end{argminie}
	where $\X$ is a complete separable metric space, $\E$ is an affine subspace of the space of continuous functions $\C^0(\X,\R^P)$, and, for every $x\in \X$, $\b F(x)$ is a non-empty closed convex subset of $\R^P$. The pointwise constraint \eqref{eq:cons_gen_F} accounts for convex constraints over the value of $\b f(x)$ whereas the global constraint \eqref{eq:cons_gen_E} can for instance force $\b f$ to be harmonic, group-invariant, and more generally takes into account linear equality constraints. We refer to \Cref{sec:examples} for examples of problem \eqref{opt-cons_gen}. The space $\C^0(\X,\R^P)$ being hard to manipulate, we approximate it by a vector-valued reproducing kernel Hilbert space (vRKHS). This choice is motivated by our main constraints of interest, \eqref{eq:cons_gen_F}, which are defined pointwise. RKHSs are precisely the Hilbert spaces of functions for which the Hilbertian topology is tighter than pointwise convergence, so that we can ensure the satisfaction of the constraints when taking limits. This is a key property of such spaces for learning tasks as recalled in \citet[Section 3]{canu2009functional}.\\
	
	For applications, we shall consider problems of the form 
	\begin{minie}|s|
		{\b f \in \HK}{\Lcal(\b f)}{\label{opt-cons_intro}}{}%
		\addConstraint{\b c_i(x)^\top\b f(x)+d_i(x)}{\ge 0,\, \forall\, x\in \Kcons_i, \, \forall \, i\in [I]:=\{1,\dots,I\}\label{eq:cons_intro_F}}
		\addConstraint{\b f}{\in \E, \label{eq:cons_intro_E}}
	\end{minie}%
	where the sets $(\Kcons_i)_{i\in [I]}$ are compact subsets of $\X$, $\b c_i\in \C^0(\X,\R^P)$, $d_i\in \C^0(\X,\R)$, and the search space $\C^0(\X,\R^P)$ is restricted to a vRKHS $\HK\subset \C^{s_K}(\X,\R^P)$ ($s_K\in\NN$), denoting by $\C^{s_K}$ the space of $s_K$-smooth functions if $\X$ is a Euclidean space. Leveraging the compactness of $\Kcons_i$, we shall sample elements within, through finite subsets $\{\tilde{x}_{i,m}\}_{m\in[M_i]}\subset \Kcons_i$. Similarly if the description of $\E$ requires an infinite number of equality constraints, we can either hard encode them in the choice of kernel $K$ (as done for translation and rotation invariant kernels in \citet{micheli_matrix-valued_2014}) or subsample the collection of equalities to define an approximating affine vector space $\hat{\E}$. Moreover $\Lcal$ could be replaced by a surrogate function $\hat \Lcal$, as when moving from population to empirical risk minimization. 
 
    Concerning the pointwise constraints, the general approach of kernel Sum-Of-Squares (kSoS) is to replace each inequality constraint by an equality constraint. This is done by adding an extra degree of freedom based on a second real-valued RKHS $\Hphi\subset \C^{s_\phi}(\X,\R)$ ($s_\phi\in\NN$) and its feature map $\phi: x\mapsto \phi(x)\in \Hphi$, considering the problem: 
	\begin{minie}|s|
		{\substack{\b f \in \HK,\\ (A_i)_{i\in[I]}\in S^+(\Hphi)^I}}{\Lcal(\b f) + \latr\sum_{i\in[I]} \tr(A_i)}{\label{opt-cons_introA}}{}%
		\addConstraint{\b c_i(\tilde{x}_{i,m})^\top\b f(\tilde{x}_{i,m})+d_i(\tilde{x}_{i,m})}{= \scalh{\phi(\tilde{x}_{i,m})}{A_i \phi(\tilde{x}_{i,m})},\, \forall\, m\in [M_i], \, \forall \, i\in [I]\label{eq:cons_introA_F}}
		\addConstraint{\b f}{\in \hat{\E}, \label{eq:cons_introA_E}}
	\end{minie}%
	where $\latr\ge 0$, and $S^+(\Hphi)$ denotes the cone of positive semidefinite operators $A$ over $\Hphi$ and $\tr(A)$ their trace. For $s_\phi\ge 2$, the extra smoothness enables the application of scattering inequalities \citep[Chapter 11]{wendland2005scattered} which improve the statistical bounds over the approximation of \eqref{opt-cons_intro} by \eqref{opt-cons_introA}. More precisely, define the fill distance as $h_{\hat{X}}:=\sup_{x\in \Kcons} \min_{m\in[M]} \|x-x_m\|$, describing how tightly $\Kcons$ is covered for $I=1$. Our main result can be informally stated as follows:
    \begin{thm}[Informal main result] 
        Under generic assumptions on the smoothness $s\in\NN$ of the functions, non-degeneracy of the constraints and a well-behaved objective function, we have that the error between the optimal value $\Lcal(\bar{\b f}^{AFF}_{\eqref{opt-cons_intro}})$ of \eqref{opt-cons_intro} and the optimal value $\Lcal(\bar{\b f}^{SDP}_{\eqref{opt-cons_introA}})$ of its approximation \eqref{opt-cons_introA} is upper bounded:
        \begin{equation*}
            |\Lcal(\bar{\b f}^{AFF}_{\eqref{opt-cons_intro}}) - \Lcal(\bar{\b f}^{SDP}_{\eqref{opt-cons_introA}})| \le C \min(h_{\hat{X}},\delta_s h_{\hat{X}}^s)
        \end{equation*}
        where $C>0$ is an explicit constant, and $\delta_s\in \{1,+\infty\}$, in particular $\delta_s=1$ if there exists $(A_i)_{i\in[I]}\in S^+(\Hphi)^I$ such that $\b c_i(x)^\top \bar{\b f}^{AFF}_{\eqref{opt-cons_intro}}(x)+d_i(x)= \scalh{\phi(x)}{A_i \phi(x)}$
    \end{thm}
    The precise statement is to be found in \Cref{thm:main_result}. In other words, our result is that sampling the constraints through kSoS as in \eqref{opt-cons_introA} always leads to a consistent approximation of \eqref{opt-cons_intro}, at worst at a slow rate, since $h_{\hat{X}}\propto M^{\nicefrac{-1}{d}}$ for uniform sampling. However, if the constraints of the optimal solution $\bar{\b f}^{AFF}_{\eqref{opt-cons_intro}}$ can be represented by a kSoS function, then this enables a fast rate of $ \mathcal{O}(M^{\nicefrac{-s}{d}})$ which mitigates the curse of dimensionality in sampling constraints.\\ 

    \tb{Comparison with previous work.} The approach and analysis we present generalizes the idea and the method proposed in \citet{rudi2020finding} for constrained convex problems resulting from the reformulation of non-convex optimization (see \Cref{sec:finding_global}), to more general kind of inequalities. Concerning the aspect of hierarchies of problems,\footnote{\label{footnote:relax-tighten}We say that problem $(\Psc_1)$ is a relaxation (resp.\ tightening) of problem $(\Psc_2)$ if they have the same objective function and the search space of $(\Psc_1)$ contains (resp.\ is contained in) that of $(\Psc_2)$. This dichotomy is paramount to our approach.} we show below in \Cref{lem:kSoS_interpolation} that, for strictly positive definite kernels $\phi$, \eqref{eq:cons_intro_F} implies the existence of $A$ satisfying \eqref{eq:cons_introA_F}. So problem \eqref{opt-cons_introA} for $\latr=0$ is a relaxation of problem \eqref{opt-cons_intro}. An alternative framework is the tightening proposed by \citet{aubin20hard,aubin2020hard_SDP}
	\begin{equation}\label{eq:cons_introSOC_F}
	\b c_i(\tilde{x}_{i,m})^\top\b f(\tilde{x}_{i,m})+d_i(\tilde{x}_{i,m}) \ge \eta_{m,i} \|\b f\|_{\HK},\, \forall\, m\in [M_i], \, \forall \, i\in [I]
	\end{equation}
	where the parameters $\eta_{m,i}>0$ only depend on the covering of $\Kcons_i$ by the discrete samples $\{\tilde{x}_{i,m}\}_{m\in[M_i]}$, on the kernel $K$ and on the constraint functions $\b c_i$ and $d_i$. Although they are in finite number, the constraints in \eqref{eq:cons_introSOC_F} are shown to be tighter than \eqref{eq:cons_intro_F}. Nevertheless the added term $\eta_{m,i} \|\b f\|_{\HK}$ guaranteeing that the constraints \eqref{eq:cons_intro_F} are satisfied also hampers the speed of convergence. The kSoS constraints \eqref{eq:cons_introA_F} leverage instead the smoothness of the functions to achieve an improved asymptotical rate as presented in \Cref{thm:main_result}. This is done at the price of not ensuring the satisfaction of the original constraints, since \eqref{opt-cons_introA} is but an approximation and not a tightening of \eqref{opt-cons_intro}. Numerically, to gain convergence speed \eqref{eq:cons_introA_F} or to achieve guarantees \eqref{eq:cons_introSOC_F} while only manipulating finitely many points, one has to pay the computational price of moving to semidefinite programming (SDP) \eqref{eq:cons_introA_F} or second-order cone programming (SOCP) \eqref{eq:cons_introSOC_F} instead of just discretizing the affine constraints \eqref{eq:cons_intro_F}.
	
	While in general we would have $\Hk \neq \Hphi$, and kSoS functions are an extra, secondary, variable, as in \eqref{opt-cons_introA}, there are situations when $f$ is only required to be nonnegative, i.e.\ \eqref{eq:cons_intro_F} corresponds to $f(x)\ge 0$. In this case, which we do not cover, kSoS functions can be used as the primary variable, $f(x)=\scalh{\phi(x)}{A \phi(x)}\ge 0$, and there is no need to discretize the constraints since kSoS functions are nonnegative. This approach was applied to estimating probability densities and then sample from the estimates \citep{marteauferey20nonparametric,marteau-ferey22sampling}.\\

	\tb{Extensions.} We can similarly consider within our framework constraints over the derivatives of $\b f$. In this case, for the derivatives to be well-defined we assume $\X$ to be a subset of $\R^d$ that is contained in the closure of its interior for the Euclidean topology. We would then replace \eqref{eq:cons_intro_F} by 
	\begin{equation}
	\b c_i(x)^\top \partial^{\b r_i} \b f(x)+d_i(x) \ge 0,\, \forall\, x\in \Kcons_i, \, \forall \, i\in [I] 
	\end{equation}
	for a multi-index $\b r_i\in\N^d$, and look for $\b f\in C^{r}(\X,\R^P)$ with $r$ the largest index appearing in the $(\b r_i)_{i\in [I]}$. Derivatives in the constraints for instance appear in the application of the kSoS approach to estimating the value function in optimal control \citep{Berthier2021InfiniteDimensionalSF}. If one wanted some components of $\b f$ to be convex, derivatives and SDP constraints would also appear, as extensively considered in \citet{aubin2020hard_SDP,muzellec21sdp}. Our approach could easily be extended to cover these cases.
    
    The article is constructed as follows. In \Cref{sec:examples}, we highlight some examples where the presented theory applies. \Cref{sec:theory} contains our main theoretical results. We illustrate these on learning vector fields in \Cref{sec:numerics}. Further technical results on closedness of sets and scattering inequalities are regrouped in \Cref{sec:closed-sets}. Finally the appendix consisting of \Cref{sec:set-valued_lsc} details more general results on selection theorems for lower semicontinuous set-valued maps to ensure that the search space of \eqref{opt-cons_gen} is non-empty.

	\section{Examples}\label{sec:examples}
	After introducing our notation and vRKHSs, we present some interesting problems (global optimization, estimation of optimal transport cost, learning vector fields) that fit our framework.
	
	\noindent\tb{Notation:}   
	Let $\N=\{0,1,\ldots\}$, $\Zp = \{1,2,\ldots\}$ and $\Rnn$ denote the set of natural numbers, positive integers and non-negative reals, respectively. We write $\iv{n_1}{n_2}=\{n_1,n_1+1,\ldots,n_2\}$ for the set of integers between $n_1, n_2\in \N$ (not to be confused with the closed interval $[a,b]$) and use the shorthand $[N]:=\iv{1}{N}$ with $N\in \N$, with the convention that $[0]$ is the empty set. The interior of a set $S\subseteq \F$ is denoted by $\mathring{S}$, its closure by $\bar{S}$, its boundary by $\partial S$. The $i^{th}$ canonical basis vector is $\b e_i$; the zero matrix is $\b 0_{d_1\times d_2} \in \R^{d_1\times d_2}$; the identity matrix is denoted by $\b I_d \in \R^{d\times d}$. The set of $d\times d$ symmetric (resp.\ positive semi-definite) matrices is denoted by $S_d$ (resp.\ $S_d^+$), $\tr(\cdot)$ denotes their trace. The closed ball in a normed space $\F$ with center $c\in \F$ and radius $r>0$ is $\BB_{\F}(c,r)=\left\{f \in \F\,:\, \left\|c - f\right\|_{\F}\le r\right\}$. Given a multi-index $\b r \in \N^d$ let $|\b r| = \sum_{j\in [d]} r_j$ be its length, and let the $\b r^{th}$ order partial derivative of a function $f$ be denoted by $\p^{\b r}f (\b x) = \frac{\p^{|\b r|}f(\b x)}{\p x_1^{r_1}\ldots \p x_d^{r_d}}$. Similarly for multi-indices $\b r,\b q \in \N^d$, let $\p^{\b r,\b q}f (\b x,\b y) = \frac{\p^{|\b r|,|\b q|}f(\b x,\b y)}{\p x_1^{r_1}\ldots \p x_d^{r_d}\p y_1^{q_1}\ldots \p y_d^{q_d}}$. For a fixed $s \in \N$, let the set of $\R^{d}$-valued functions on $\X$ with continuous derivatives up to order $s$  be denoted by $\C^s\left(\X,\R^{d}\right)$. 
	The set of $\R^{d_1\times d_2}$-valued functions on $\X\times \X$ for which  $\p^{\b r,\b r}f$ exists and is continuous up to order $|\b r|\le s\in \N$  is denoted by $\C^{s,s}\left(\X\times \X,\R^{d_1\times d_2}\right)$. The domain $\Dom(f)$ of a function is the set of points where it has finite values.%

	\tb{vRKHS:}	A function $K: \X \times \X \rightarrow \R^{Q \times Q}$ is called a positive semidefinite matrix-valued kernel on $\X$ if $K(\b x,\b x')=K(\b x',\b x)^\top$ for all $\b x,\b x' \in \X$ and $\sum_{i,j\in[N]} \b v_i^\top K(\b x_i,\b x_j)\b v_j\ge 0$ for all  $N\in \Zp$, $\left\{\b x_n\right\}_{n\in [N]} \subset \X$	and $\{\b v_n\}_{n\in [N]} \subset \R^Q$.  For $\b x\in \X$, let $K(\cdot,\b x)$ be the mapping $\b x' \in \X \mapsto K\left(\b x',\b x\right) \in \R^{Q\times Q}$. Let $\FK$ denote the vRKHS associated to the kernel $K$ \citep[see][and references therein]{carmeli10vector}; we use the shorthand $\left\|\cdot\right\|_K := \left\|\cdot\right\|_{\FK}$ and $\left<\cdot,\cdot\right>_K := \left<\cdot,\cdot\right>_{\FK}$ for the norm and the inner product on $\FK$. The Hilbert space $\FK$ consists of $\X \rightarrow \R^Q$ functions for which (i) $K(\cdot,\b x)\b c \in \FK$ for all $\b x\in\X$ and $\b c \in \R^Q$, and (ii) $\left<\b f,K(\cdot,\b x)\b c\right>_{K}=\left<\b f(\b x),\b c\right>$ for all $\b f\in \FK$, $\b x\in\X$ and $\b c\in \R^Q$. The first property of vRKHSs describes the basic elements of $\FK$, the second one is called the reproducing property; this property can be extended to function derivatives \citep[Lemma 1]{aubin2020hard_SDP}. %
	Constructively, $\FK = \overline{\Sp}\left\{K(\cdot,\b x)\b c\,:\, \b x \in \X, \b c \in \R^Q\right\}$ where $\Sp$ denotes the linear hull of its argument and the bar stands for closure w.r.t.\ $\|\cdot\|_K$. 
	
	\subsection{Finding global minima}\label{sec:finding_global}
	We follow here the construction in \citet{rudi2020finding} which develops a subset of the framework elaborated in this paper. The goal in this context is to solve
     $$g_* = \inf_{x \in \Omega} g(x),$$
     where $g \in \C^s(\R^d)$ with $s \ge 2$, $\Omega$ is an open bounded subset of $\R^d$, with $d \in \N$ and there exists at least a minimizer $x_* \in \Omega$. A simple convex reformulation of the problem above, is as follows
    \begin{mini*}|s|
		{c \in \R}{-c}{}{}%
		\addConstraint{g(x) - c}{\ge 0,\, \forall\, x\in \Omega}
	\end{mini*}
    The problem above gives exactly $g_*$ as solution, indeed $c=g_*$ satisfies the constraints and all $c > g_*$ would violate the constraints in a neighbourhood of $x_*$ which is in $\Omega$. The problem above is convex, since both the functional to minimize and the constraints are linear in $c$, and fits our problem of interest \eqref{opt-cons_intro}, by setting $k(x,y)=1$ encoding the constant functions, $I=1$, $c(x)=-1$, $d(x)=g(x)$. By applying the construction we propose in this paper, we obtain the following convex approximation of the original non-convex minimization problem
	\begin{mini*}|s|
		{\substack{c \in \R,\\ A \in S^+(\hh_\phi)}}{-c + \latr \tr(A)}{}{}%
		\addConstraint{g(\tilde{x}_{m}) - c}{= \scalidx{\phi(\tilde{x}_{m})}{A \phi(\tilde{x}_{m})}{\hh_\phi},\, \forall\, m\in [M]}
	\end{mini*}
where $\tilde{x}_1, \dots, \tilde{x}_M$ are a list of points chosen appropriately (e.g.\ uniformly at random), to cover $\Omega$ and $\hh_\phi$ is a rich enough RKHS (as a Sobolev space). Under some assumptions on $g$, such as isolated minima with positive Hessians, and on $\hh_\phi$, the solution of the problem above converges to $g_*$ as $M$ goes to infinity and with convergence rates that are almost optimal, as discussed in \citet{rudi2020finding}.%

Our general analysis covers the algorithm above as a particular case, and allows to derive slightly better rates. Indeed, solving the global optimization problem $g_* = \inf_x g(x)$ up to error $\eps$ is very expensive. Typical algorithms to solve global optimization have computational complexity of the order $O(\eps^{-d/2})$ \citep[see Section 1][for a discussion]{rudi2020finding}. However, the known computational lower-bound is in the order of $O(\eps^{-d/s})$ \citep{Novak1988}, breaking the curse of dimensionality when $g$ is very regular, i.e. $s \gg 1$. For example, if $s \geq d$, it means that an optimal algorithm would take time at most $O(\eps^{-1})$, in other words to double the precision it would just need double the time. \citet{rudi2020finding} show that, by applying the technique under discussion, the algorithm has computational complexity $O(\eps^{-3.5 d/(s- d/2-2)})$. While this rate is not optimal (mainly for the factor $3.5$ in the exponent), again it has the same qualitative effect of breaking the curse of dimensionality, indeed for $g$ very regular, e.g.\ $s \geq 4d + 2$ the algorithm has a rate of $O(\eps^{-1})$.
	
	\subsection{Optimal transport with kSoS}
Another problem that fits perfectly the framework under study is the one of estimating the optimal transport cost $W_c(\mu, \nu)$  between two smooth probability densities $\mu, \nu$ defined over a bounded set $\Omega \subset \R^d$ with $d \in \N$. The Kantorovich dual formulation reads as follows \citep[see Chapter 5][for an introduction]{Villani2009}
\begin{mini*}|s|
		{u, v \in C(\R^d)}{ -\int u(x) d\mu(x) - \int v(y) d\nu(y)}{}{W_c(\mu, \nu) = }%
		\addConstraint{u(x) + v(y)}{ \leq c(x,y),\, \forall\, x, y \in \Omega}
	\end{mini*}%
where $c(x,y)$ is the cost of transporting a unit of mass from $x$ to $y$ and is assumed to be a smooth function, e.g.\ $c(x,y) = \|x-y\|^2$ (in general $x$ and $y$ can be taken in different sets with few modifications). In this case, applying the proposed method to $\HK=\hh_k\times \hh_k$, $K(x,y)=k(x,y) \Id_2$ would result in 
\begin{mini*}|s|
		{\substack{u, v \in \hh_k, \\ A \in S^+(\hh_\phi)}}{ -\int u(x) d\mu(x) - \int v(y) d\nu(y) ~+~ \latr \tr(A)} {}{}%
		\addConstraint{u(\tilde{x}_m) + v(\tilde{y}_m) - c(\tilde{x}_m, \tilde{y}_m)}{ = \scalh{\phi(\tilde{x}_{m}, \tilde{y}_{m})}{A \phi(\tilde{x}_{m}, \tilde{y}_{m})},\, \forall\, m \in [M]}
\end{mini*}
 where $\Hphi$ is a RKHS on $\Omega \times \Omega$ and $\hh_k$ is a RKHS on $\Omega$. Moreover, if we cannot perform the integral in $\mu, \nu$, but we have some samples $x_1,\dots, x_n$ and $y_1,\dots, y_n$, that are sampled i.i.d.\ from $\mu$ and $\nu$, we can also approximate the integral as $\int u(x) d\mu(x) \approx \frac{1}{n} \sum_{i} u(x_i)$, obtaining the following problem,
 \begin{mini*}|s|
		{\substack{u, v \in \hh_k, \\ A \in S^+(\hh_\phi)}}{ -\frac{1}{n} \sum_{i=1}^n u(x_i) - \frac{1}{n} \sum_{i=1}^n v(y_i) ~+~ \latr \tr(A)}{}{ }%
		\addConstraint{u(\tilde{x}_m) + v(\tilde{y}_m) - c(\tilde{x}_m, \tilde{y}_m)}{ = \scalh{\phi(\tilde{x}_{m}, \tilde{y}_{m})}{A \phi(\tilde{x}_{m}, \tilde{y}_{m})},\, \forall\, m \in [M]}
	\end{mini*}
 The solution of the problem above converges to $W_c(\mu, \nu)$ as $M$ and $n$ go to infinity, as discussed in \citet{vacher2021ot}, that first presented this application of kSoS techniques. Our general analysis covers the algorithm above as a particular case, and provides the same optimal rates, but with explicit constants. Moreover, as in the non-convex optimization case of \Cref{sec:finding_global}, computing the Wasserstein distance is also in general a difficult problem suffering from the curse of dimensionality. However, the paper \citep{vacher2021ot} showed that the approach under study achieves an error $\eps$ with a computational complexity that breaks the curse of dimensionality in the exponent, depending on the degree of regularity of the problem (in this case, $\mu, \nu$ should have $s$-times smooth densities w.r.t.\ the Lebesgue measure). When $s \gg d$, the resulting complexity is then dimension-independent, although the constants in the rate may depend exponentially on $d$.
	
	\subsection{Learning dynamical systems with shape constraints}	
	
	In this novel application of the framework, our goal is to learn a continuous-time autonomous dynamical system	of the form
	\begin{equation}\label{eq:dyn_sys}
	\dot{x}(t)=\b f(x(t))
	\end{equation}
	over a set $\X\subseteq \R^d$, based on a few noisy observations of some trajectories of the system \eqref{eq:dyn_sys}, forming a training set
	\begin{equation}\label{eq:training_set}
	\D:=\{(x_n,\b y_n), \, n\in[N]\}\subset(\X\times\R^{d})^N.
	\end{equation}
	If the velocity vectors $\b y_n$ are not directly available from measurements, they can be estimated based on the observed trajectories. %
	However we expect that the modeler has access to some extra information on the properties of the desired $\b f$ beyond the samples $\D$. The vector field $\b f$ may for instance be known to leave invariant some subsets of $\R^d$, or to be the gradient of some unknown potential function. For controlled systems, $\b f$ as in \eqref{eq:dyn_sys} can correspond to a given closed-loop control, i.e.\ for a given $\tilde{\b f}:\X\times \U\rightarrow\R^d$, with a control set $\U\subset\R^m$, there exists a function $\bar{\b u}:\X \rightarrow \U$ such that
	\begin{equation}\label{eq:control}
	\b f(x)=\tilde{\b f}(x,\bar{\b u}(x)), \, \forall \, x\in \X.
	\end{equation}
	We thus look for an estimation procedure of a vector field $\hat{\b f}$ that is both compatible with the observations, in the sense that $\b y_n \approx \hat{\b f}(x_n)$, and with the qualitative priors imposed on its form, such as \eqref{eq:control}. These priors are a form of side information that is incorporated to improve both the estimation procedure of $\hat{\b f}$ and its interpretability in the context of dynamical systems \citep{ahmadi23learning}. In statistics, they are known as shape constraints, such as assuming non-negativity, monotonicity or convexity of some components of the estimation \citep{guntuboyina18nonparametric}. Many of these constraints can be written as requiring that $\hat{\b f}(x)$ belongs to a set $\b F_0(x)\subseteq \R^d$,\footnote{We then say that $\hat{\b f}$ is a \emph{selection} of the set-valued map $\b R:\X \leadsto \R^d$.}%
	\begin{equation}\label{eq:constr_set-valued}
	\hat{\b f}(x)\in \b F_0(x), \, \forall \, x\in \X.
	\end{equation}%
	For instance, if we want to interpolate some points, we can take $\b F_0(x_i)=\{\b y_i\}$; if we want the first component $\hat{f}_1$ to be nonnegative over a set $\Kcons\subset \R^d$, then we can pose $\b F_0(x)=\{\b y\,|\, y_1\ge 0\}$ for all $x\in\Kcons$ and $\b F_0(x)=\R^d$ elsewhere; for controlled systems as in \eqref{eq:control}, $\b F_0(x)=\tilde{\b f}(x,\U)$. Constraints of the form \eqref{eq:constr_set-valued} correspond to local pointwise requirements of $\hat{\b f}$. We could also consider nonlocal constraints such as parity of $\hat{\b f}$, i.e.\ $\hat{\b f}(-x)=\hat{\b f}(x)$, or global requirements such as $\hat{\b f}=\nabla V$ for some $V:\X\rightarrow \R$. Each of these side information corresponds to constraining $\hat{\b f}$ to belong to a set of functions $\E$. 
	
	Moreover we require $\hat{\b f}$ to belong to a set of functions $\b f$ such as the spaces of Lipschitz or of $\C^1$-smooth functions over $\X$ in order to perform some minimal, Cauchy-Lipschitz style, analysis of \eqref{eq:dyn_sys}.\footnote{Given some trajectories, there does not always exist an autonomous system such as \eqref{eq:dyn_sys} if time is involved. For instance a railroad switch may change with time bringing trains to different directions even though they pass through the same point. Non-autonomous systems can be seen as a succession of autonomous systems on very short time intervals and would require many observations at every given time to be learnt.} \citet{ahmadi23learning} proposed to use polynomial sum-of-squares to fit the constraints and polynomial functions for $\b f$. However polynomials of degree two and above induce finite-time explosion of trajectories. Moreover the optimization problem becomes quickly computationally intractable when the degree of the polynomials increases, due to the combinatorial number of monomials. Instead kernels are much smoother, and often chosen so as to vanish at infinity. Applying matrix-valued kernels to recover vector fields has been already considered \citep[see e.g.][and references within]{Baldassarre2012multi}, but, up to our knowledge, not with infinitely many constraints to the exception of \citet{aubin2020hard_SDP}. To assess the performance of the estimator, one can consider the $L^2$ or $L^\infty$ error between vector fields or an angular error, the distance between reconstructed trajectories, the amount of violation of constraints,\dots We present a numerical application on estimating vector fields in \Cref{sec:numerics}.

	\section{Theoretical results}\label{sec:theory}%
 	
	In what follows, to simplify the presentation, we take $\hat \E=\E=\HK$, $\hat \Lcal = \Lcal$ in \eqref{opt-cons_intro}-\eqref{opt-cons_intro} and consider a single constraint set $\Kcons_i=\Kcons$ ($\forall i\in[I]$), focusing only on the approximation of the inequality constraints. We shall also consider the following assumptions, which essentially require some smoothness of the functions, non-degenerate constraints and a well-behaved objective function.
	\newcounter{contlist}
	\begin{assumplist}
		\item \label{ass:continuous_kernel} The set $\X$ is contained in $\R^d$ and in the closure of its interior. The kernels satisfy $K(\cdot,\cdot) \in \C^{s_K,s_K}\left(\X\times \X,\R^{P\times P}\right)$, $k_\phi(\cdot,\cdot) \in \C^{s_\phi,s_\phi}\left(\X\times \X,\R\right)$ for some $s_K,s_\phi \in \N$.%
		
		\item \label{ass:continuous_cons} The constraint functions satisfy $\b C(\cdot) \in \C^{s}\left(\X,\R^{I\times P}\right)$ and  $\b d(\cdot) \in \C^{s}\left(\X,\R^{I}\right)$ for some $s\in \N$.
		
		\item \label{ass:constraint_set} The set $\Kcons$ is a compact set with samples $\hat{X}=\{x_m\}_{m\in[M]}\subset \Kcons$ chosen such that $\b K_{\phi,\hat{X}}=[k_\phi(x_i,x_j)]_{i,j\in[M]}$ is invertible.
		
		\item \label{ass:nondegen_constraint} For any $x\in\Kcons$ there exists $\b z_x\in\R^P$ such that $\b C(x)\b z_x>0$ and the functions of $\HK$ when restricted to $\Kcons$ are dense in $C^0(\Kcons,\R^P)$.
		
		\item \label{ass:smooth_objective} The objective $\Lcal$ is lower bounded over $\HK$, has full domain, i.e.\ $\HK \subset \Dom(\Lcal)$, and its sublevel sets are weakly compact\footnote{A subset of a Hilbert space is weakly compact if every sequence in it has a weakly converging subsequence to an element of the set \citep[see][Section 2.4.4]{attouch14variational}.} in $\HK$. 
		
	\end{assumplist}
 \begin{assumplistp}[start =4]
            \item\label{ass:nondegen_constraint_short}  \rv{ There exists $g\in\HK$ such that, for any $x\in\Kcons$, we have $\b C(x)g(x)\ge \b 1$.}
		
	\end{assumplistp}
	\noindent \tb{Remarks on the assumptions.} If no derivatives are considered, i.e.\ $s=0$, then $\X$ in \Cref{ass:continuous_kernel} can be any complete separable metric space, and bounds on derivatives can be replaced by Lipschitz bounds, provided the functions under study are Lipschitz on $\Kcons$. Smoothness as in \Cref{ass:continuous_kernel,ass:continuous_cons} is required for lower semicontinuity of the objective and to trigger the fast approximation behaviour of the kSoS estimate. Compactness in \Cref{ass:constraint_set} allows to consider the filling distance and to ensure that sampling eventually covers the whole set. Any strictly positive definite (resp.\ universal) kernel satisfies the second part of \Cref{ass:constraint_set} (resp.\ \Cref{ass:nondegen_constraint}). The Gaussian kernel is for instance both strictly positive definite and universal. Our assumptions are more general and lighter than those of \citet{rudi2020finding} which essentially restricted the kernel $\phi$ to be of the Sobolev type, as their proofs require $\Hphi$ to be a multiplicative algebra. \rv{On the other hand, \Cref{ass:nondegen_constraint_short} is sufficient to have a non-empty constraint set but it may be hard to show except in simple cases, such as the constant functions of \Cref{sec:finding_global} which in turn do not satisfy \Cref{ass:nondegen_constraint} since the density is lacking. In \Cref{thm:well-posed_cons} below, we show that \Cref{ass:continuous_kernel,ass:continuous_cons,ass:nondegen_constraint} together imply \Cref{ass:nondegen_constraint_short}.} Finally, \Cref{ass:smooth_objective} is a classical assumption in infinite-dimensional optimization to ensure existence of solutions, as it is satisfied by convex lower semicontinuous coercive functions \citep[Definition 3.2.4]{attouch14variational}.\\
	
	We shall manipulate various forms of inequality constraints. We thus write our generic problem \eqref{opt-LQR_gen}, with objective $\Lcal$ and constraints defined through the inequality constraint set $\V^*$, as follows
	\begin{argmini}
		{\substack{\b f(\cdot)\in \HK}}{\Lcal(\b f(\cdot))}{\tag{$\Psc_{*}$}}{ \bar{\b f}^* (\cdot)\in\label{opt-LQR_gen}}
		\addConstraint{\b f(\cdot)}{\in \V^*.}
	\end{argmini}
	We assumed that all the equality constraints \eqref{eq:cons_gen_E} on $\b f$ were linear, and moved to the definition of $\HK$. The inequality constraints \eqref{eq:cons_gen_F} should remain convex even when turned into equalities, so we will stick to the case of affine constraints \eqref{eq:cons_intro_F}. Unlike in \eqref{opt-cons_introA}, we will constrain $\tr(A)$ rather than penalize it. Both problems are of course equivalent, but constraints are more amenable for theoretical results, while penalization is more practical computation-wise.
	
	For any $\epsArr\in\R^I$, $M_f\ge 0$, $M_A\ge 0$, given $\hat{X}\subset \Kcons\subset \X$ as in \Cref{ass:constraint_set}, we shall consider the following constraints
	\begin{align}
	\V^{AFF}_0&:=\{\b f(\cdot)\in\HK \,|\,\b C(x) \b f (x)+ \b d(x)\ge  \b 0,\,\forall \, x\in \Kcons\}\label{eq:def_aff_constraints}\\
	\V^{SDP}_{\epsilon}&:=\{\b f(\cdot)\in\HK \,|\, \exists (A_i)_{i\in[I]}\in S^+(\Hphi)^I,\nonumber\\ &\hspace{3cm}\b C (x_m) \b f (x_m)+ \b d(x_m)=\epsArr+[\langle \phi(x_m),  A_i \phi(x_m)\rangle_{\Hphi}]_{i\in[I]},\, \forall \, m\in [M] \},\label{eq:def_sdp_constraints} \\
	\V^{SDP}_{\epsilon,rest}&:=\{\b f(\cdot)\in \V^{SDP}_{\epsilon} \text{ with }  \max_{i\in[I]}\tr(A_i)\le M_A,\, \|\b f\|_K\le M_f \} \label{eq:def_sdp_constraints_rest}.
	\end{align}
	We also introduce three quantities which will be useful to study the properties of the SDP approximation \eqref{opt-cons_introA}. They are defined for $s\in\NN$ such that $s\le s_\phi$, a collection of points $\hat{X}=\{x_m\}_{m\in[M]} \subset \X$, a bounded set $\Omega$, and any function $g\in C^s(\X,\R)$, as follows
	\begin{align}
	&\text{(fill distance)}\quad &&h_{\hat{X}, \Omega}:=\sup_{x\in \Omega} \min_{m\in[M]} \|x-x_m\|; \label{eq:fill_distance} \\
	&\text{($C^s$ seminorm)}\quad &&|g|_{\Omega,s}:=\sup_{x\in \Omega} \max_{\bm \alpha \in \N^d,\, |\bm\alpha|=s} |\p^{\bm\alpha} g(x)|; \label{eq:deriv_bound}\\
	&\text{(bound on kernel derivatives)}\quad &&D_{\Omega,s}^2:= \sup_{x,y\in\Omega} \max_{\bm \alpha \in \N^d,\, |\bm \alpha| \le s} |\p_{x}^{\bm\alpha}\p_{y}^{\bm\alpha}k_\phi(x,y)|. \label{eq:kernel_deriv_bound}
	\end{align}
	\tb{Roadmap of the results.} We first prove in \Cref{thm:well-posed_cons} that our assumptions ensure the well-posedness, i.e.\ existence of solution, of each of our three problems. We then show in \Cref{thm:main_result} that, for smooth positive definite kernels, the SDP approximation ($\Psc_{0}^{SDP}$) is a relaxation whereas, for well-chosen $(M_f,M_a,\epsilon)$, the SDP approximation ($\Psc_{\epsilon,rest}$) is a tightening. Exhibiting a nested sequence of problems allows to derive rates of convergence. Moreover the amount of perturbation of the constraints proves to be crucial to bound easily the approximation error. More precisely, the result concerning the relaxation is based on \Cref{lem:kSoS_interpolation} which can be seem as a converse of \citet[][Lemma 3, p15]{rudi2020finding}. The weak closedness of the sets, which is used in proving the well-posedness, is stated in the more technical \Cref{lem:closed_sdp_constraints} in \Cref{sec:closed-sets}. \Cref{sec:set-valued_lsc} contains two further results on smooth selections in constraint sets, \Cref{lem:lsc_affine constraints} and \Cref{cor:selection_affine_cons}, which ensure that the constraint sets are not empty.\footnote{More generally, returning to the original problem \eqref{opt-cons_gen} of which \eqref{opt-cons_intro} is an instance, well-posedness is guaranteed for $\E=\C^0(\X,\R^P)$, as detailed in \Cref{sec:set-valued_lsc}, provided that i) the set-valued constraint map $\b F:\X\leadsto \R^P$ is lower semicontinuous with closed, convex values, that ii) the objective $\Lcal$ is lower semicontinuous and lower bounded on the considered set of constraints.}

	\begin{lemma}[kSoS interpolation]\label{lem:kSoS_interpolation} Given $(x_m,y_m)_{m\in[M]}\in (\X \times \R_+)^M$ s.t.\ $\b K_{\phi,\hat{X}}=[k_\phi(x_i,x_j)]_{i,j\in[M]}$ is invertible, then there exists $A\in S^+(\Hphi) $ such that $y_m=\langle \phi(x_m),  A \phi(x_m)\rangle_{\Hphi}$, and $\tr(A)=\tr(\diag(\b y)\b K_{\phi,\hat{X}}^{-1})\le \min(\|\b y\|_\infty \tr(\b K_{\phi,\hat{X}}^{-1}), \|\b y\|_1 \lambda_{min}(\b K_{\phi,\hat{X}}))$ where $\lambda_{min}(\b K_{\phi,\hat{X}})$ is the smallest eigenvalue of $\b K_{\phi,\hat{X}}$.
	\end{lemma}
	\begin{proof}
		Let $\widehat{S}:\Hphi \rightarrow \RR^M$ be the linear operator acting as follows
		$$
		\widehat{S} g=\left(\left\langle\phi(x_{1}), g\right\rangle_{\Hphi}, \ldots,\left\langle\phi(x_{M}), g\right\rangle_{\Hphi}\right) \in \RR^M, \quad \forall g \in \Hphi
		$$
		Consider the adjoint of $\widehat{S}$, $\widehat{S}^{*}: \RR^M \rightarrow \Hphi$, defined consequently as $\widehat{S}^{*} \beta=\sum_{m\in[M]} \beta_m \phi\left(x_m\right)$ for $\bm \beta \in \RR^M$. Note, in particular, that $\b K_{\phi,\hat{X}}=\widehat{S} \widehat{S}^{*}$ and that $\widehat{S}^{*} e_m=\phi\left(x_m\right)$, where $e_m$ is the $m$-th element of the canonical basis of $\RR^M$. Fix the continuous linear operator $A:=\widehat{S}^{*} \b K_{\phi,\hat{X}}^{-1} \diag(\b y) \b K_{\phi,\hat{X}}^{-1} \widehat{S}$ acting from $\Hphi$ to $\Hphi$. The operator $A$ is clearly self-adjoint as the kernel is positive semidefinite, hence $\b K_{\phi,\hat{X}} \in S^+(\RR^M)$. Moreover, for any $g\in\Hphi$, $\left\langle g, A g\right\rangle_{\Hphi}=\left\langle \b V , \diag(\b y) \b V\right\rangle_{\RR^M}$ with $\b V:=\b K_{\phi,\hat{X}}^{-1} \widehat{S} g$ and $\b y=(y_m)_m\ge 0$. So $A$ is positive semidefinite. In particular,
		\begin{equation*}
		\langle \phi(x_m),  A \phi(x_m)\rangle_{\Hphi}=\langle \widehat{S}^{*} e_m,  A \widehat{S}^{*} e_m\rangle_{\Hphi}=\langle e_m,\widehat{S}  A \widehat{S}^{*} e_m\rangle_{\RR^M}=\langle e_m,\diag(\b y) e_m\rangle_{\RR^M}=y_m.
		\end{equation*}
		
		Using the cyclicity of the trace and that $\b K_{\phi,\hat{X}}=\widehat{S} \widehat{S}^{*}$, we obtain that $\tr(A)=\tr(\widehat{S}\widehat{S}^{*} \b K_{\phi,\hat{X}}^{-1} \diag(\b y) \b K_{\phi,\hat{X}}^{-1})=\tr(\diag(\b y)\b K_{\phi,\hat{X}}^{-1})$. To prove that $\tr(A)\le \min(\|\b y\|_\infty \tr(\b K_{\phi,\hat{X}}^{-1}), \|\b y\|_1 \lambda_{min}(\b K_{\phi,\hat{X}}))$, we just use that for any matrices $B, C\in S^+(\RR^M)$, $\tr(BC)\le \lambda_{max}(B) \tr(C)$ \citep[e.g.][Lemma 1]{sheng1986trace}.
	\end{proof}

	\begin{theorem}[Well-posedness of the problems]\label{thm:well-posed_cons} We have the following:
		
		\begin{enumerate}[label={\arabic*)},labelindent=0em,leftmargin=1em,topsep=0.2cm,partopsep=0cm,parsep=0cm,itemsep=2mm]
			
			\item\label{it_aff_subset_eps} Provided \Cref{ass:constraint_set} is satisfied, we have the inclusion $\V^{AFF}_0\subset \V^{SDP}_0$ and we can fix $M_f(\bar{\b f}^{AFF}_0) \in \R_+$ and some $M_A(\bar{\b f}^{AFF}_0,\hat{X})\in \R_+$ such that, if $\bar{\b f}^{AFF}_0$ exists, then $\bar{\b f}^{AFF}_0\in \V^{SDP}_{0,rest}$.
			
			\item\label{it_nonempty}  If \Cref{ass:continuous_kernel,ass:continuous_cons,ass:nondegen_constraint}  hold, then \Cref{ass:nondegen_constraint_short} holds and we can fix $\b g_1= \argmin{\b g \in H_1}\|\b g\|_{\HK}$ where $H_1:=\{\b g\in\HK \, |\, \b C(x) \b g(x) \ge \b 1,\, \forall x\in\Kcons\}$. Moreover the set $\V^{AFF}_0$ is then non-empty and, if \Cref{ass:constraint_set} also holds, so is $\V^{SDP}_0$.
			
			\item\label{it_solutions_existence} Under \Cref{ass:smooth_objective,ass:continuous_cons,ass:nondegen_constraint,ass:constraint_set}, $\bar{\b f}^{AFF}_0$ exists and, for $M_f(\bar{\b f}^{AFF}_0)$ and $M_A(\bar{\b f}^{AFF}_0,\hat{X})$ as in \ref{it_aff_subset_eps}, $\bar{\b f}^{SDP}_{0,rest}$ also exists.
		
		\end{enumerate}
		
	\end{theorem}
	\begin{proof}
		
		\tb{Proof of \ref{it_aff_subset_eps}:} If $\b K_{\phi,\hat{X}}=[k_\phi(x_i,x_j)]_{i,j\in[M]}$ is invertible, by \Cref{lem:kSoS_interpolation}, for every given $\b f \in \V^{AFF}_0$ and $i\in[I]$, there exists an operator $A_i\in S^+(\Hphi)$ interpolating the values $y_{m,i}=\b c_i(x_m)^\top \b f (x_m)+ d_i(x_m)\ge 0 $ as a kernel Sum-of-Squares. Thus $\V^{AFF}_0\subset \V^{SDP}_0$. In particular, if it exists, $\bar{\b f}^{AFF}_0\in \V^{AFF}_0\subset \V^{SDP}_0$. Thus we just have to set $M_f=\|\bar{\b f}^{AFF}_0\|_K$ and $M_A(\bar{\b f}^{AFF}_0,\hat{X})= \tr(\b K^{-1}) \max_{m\in[M], i\in[I]} |y_{m,i}|$, as discussed in \Cref{lem:kSoS_interpolation}, in order to have $\bar{\b f}^{AFF}_0\in \V^{SDP}_{0,rest} $.\\

		\tb{Proof of \ref{it_nonempty}:}\footnote{In \citet[][Appendix]{aubin2020hard_control}, the space $\HK$ was given by linearly controlled functions. Hence some more specific arguments had to be given, i.e.\ inward-pointing conditions, since one could not leverage generic density properties of the function space.} \Cref{cor:selection_affine_cons} written with $\zeta=1$ provides a function $\b f_1\in C^0 (\X,\R^P)$ such that $\b f_1(x)+ \BB_{\R^P}(0,1) \subset \{\b y\,|\, \b C(x) \b y \ge \b 1\}$ for all $x\in\Kcons$. The latter is equivalent to saying that $\b 1 \le \inf_{\b u\in \BB_{\R^P}(0,1)} \b C(x) (\b f_1(x)+\b u)$. Since in \Cref{ass:nondegen_constraint}, we assumed the functions of $\HK$ when restricted to $\Kcons$ to be dense in $C^0(\Kcons,\R^P)$, we can find $\b g\in\HK$ such that $\sup_{x\in \Kcons}\|\b f_1(x)-\b g(x)\|_{\infty}\le \frac{1}{\sqrt{P}}$. Consequently, for any $x\in\Kcons$, $(\b f_1-\b g)(x)\in \BB_{\R^P}(0,1)$, i.e.\ it is an element over the set over which the infimum is taken, so
		\begin{align*}
        \b C(x) \b g(x)&=\b C(x) (\b f_1(x)+\b g(x)-\b f_1(x))\ge \inf_{\b u\in \BB_{\R^P}(0,1)} \b C(x) (\b f_1(x)+\b u)\ge \b 1,
		\end{align*}
		hence $\b g\in H_1$ which is thus non-empty. Since $K$ is a reproducing kernel, for all $x\in\X$, $H_{1,x}:=\{\b g\in\HK \, |\, \b C(x) \b g(x) \ge \b 1\}$ is a closed convex subset of $\HK$. So is $H_1=\cap_{x\in\Kcons}H_{1,x}$ as an intersection of such sets. The element $\b g_1$ thus exists as it is the projection of $\b 0$ onto $H_1$.
		
		For any $R>0$ define $\b g_R =R \b g_1\in\HK$. For $R=\sup_{x\in \Kcons}\|\b d(x)\|_\infty$, 
		\begin{align*}
		\b C(x) \b g_R(x)&\ge \sup_{x\in \Kcons}\|\b d(x)\|_\infty \b 1 \ge -\b d(x),
		\end{align*}
		so $\b g_R\in\V^{AFF}_0$ which is thus non-empty. If \Cref{ass:constraint_set} also holds, as shown in \ref{it_aff_subset_eps}, we obtain that $\V^{SDP}_0\supset \V^{AFF}_0\neq \emptyset $.\\

		\tb{Proof of \ref{it_solutions_existence}:} Under \Cref{ass:smooth_objective,ass:continuous_cons,ass:nondegen_constraint,ass:constraint_set}, by \ref{it_aff_subset_eps}, we have that $\emptyset \neq \V^{AFF}_0\subset \V^{SDP}_0$ and $\V^{SDP}_{0,rest} \neq \emptyset$. Moreover \Cref{lem:closed_sdp_constraints} gives that these sets are all weakly closed in $\HK$. Denote generically by $\V^*$ any of these sets as was done in \eqref{opt-LQR_gen}. Since $\Lcal$ has full domain and is lower-bounded, it has a finite smallest value $\bar{v}^*$ over $\V^*$. By definition of the infimum, we can fix $(\b f_n)_{n\in\N}\in (\V^*)^\N$ such that  the sequence $\Lcal(\b f_n)\in[\bar{v}^*,\bar{v}^*+1] $ converges to $\bar{v}^*$. By \Cref{ass:smooth_objective}, the sublevel set \rv{$\F^*=\{\b f\in\HK \, | \, \Lcal(\b f)\le \bar{v}^*+1\}$} is weakly compact. As $\V^*$ is weakly closed again by \Cref{lem:closed_sdp_constraints}, $\F^*\cap\V^*$ is then also weakly compact, whence $(\b f_n)_{n\in\N}\in (\F^*\cap\V^*)^\N$ has a weakly converging subsequence to some $\bar{\b f}^*$ which achieves the infimum, i.e.\ $\Lcal(\bar{\b f}^*)=\bar{v}^*$. So all the infima exist as claimed.
	\end{proof}
	For a small fill distance  $h_{\hat{X}, \Omega}$, and thus a large number $M$ of sample points, we show in \Cref{thm:main_result} below that there is a nested sequence between $\V^{SDP}_{\sigma,rest}$,  $\V^{AFF}_0$ and $\V^{SDP}_0$, defined respectively in \eqref{eq:def_sdp_constraints_rest}, \eqref{eq:def_aff_constraints} and \eqref{eq:def_sdp_constraints}, and where $\sigma\in\R_+$ is related to the parameters of the problem and decreases with $h_{\hat{X}, \Omega}$. The role of $\sigma$ to derive rates of convergence was partly anticipated in \citet[][Theorem 4]{rudi2020finding}, whereas the idea of using nestedness for studying similar constraints perturbation was advocated in \citet[][Proposition 1]{aubin2020hard_control}.

	\begin{theorem}[Main result]\label{thm:main_result} Suppose \Cref{ass:continuous_cons,ass:nondegen_constraint,ass:constraint_set,ass:continuous_kernel}, \rv{or \Cref{ass:continuous_cons,ass:constraint_set,ass:continuous_kernel} and \Cref{ass:nondegen_constraint_short}}, hold, with  $1 \le s \le \min(s_K,s_\phi)$. For $M_f(\bar{\b f}^{AFF}_0)$ and $M_A(\bar{\b f}^{AFF}_0,\hat{X})$ given by \Cref{thm:well-posed_cons}-\ref{it_aff_subset_eps}, define $C_{K,s}=\max_{x\in\Kcons, i\in[I]} \max_{\bm \alpha \in \N^d,\, |\bm \alpha|= s}\|\p_x^{\bm \alpha}[K(\cdot,x)\b c_i(x)]\|_{\HK}$ and \footnote{See \eqref{eq:fill_distance}, \eqref{eq:deriv_bound} and \eqref{eq:kernel_deriv_bound} for the definition of the quantities}
		\begin{equation}\label{eq:err_SDP_1}
		\sigma_{slow}:=\left(C_{K,1} M_f(\bar{\b f}^{AFF}_0)+\max_{i\in[I]}|d_i|_{\Kcons, 1}\right) \cdot h_{\hat{X},\Kcons}.
		\end{equation}
	Moreover, for  $s\ge 2$, if $\Kcons= \cup_{x\in S}\BB(x,r)$ for a given bounded set $S\subset \X$ and some $r>0$, and, setting $\Omega:= \cup_{x\in S}\dBB(x,r)$, we have that $\hat{X}=\{x_m\}_{m\in[M]}\subset \Omega$ with $h_{\hat{X}, \Omega}\le r \min(1,\frac{1}{18(s-1)^2})$, then we define
			\begin{equation}\label{eq:err_SDP_s}
				\sigma_{fast}:=C_0\left(C_{K,s} M_f(\bar{\b f}^{AFF}_0)+\max_{i\in[I]}|d_i|_{\Kcons, s}+2^s D_{\Omega,s}^2 M_A(\bar{\b f}^{AFF}_0,\hat{X})\right) \cdot(h_{\hat{X},\Kcons})^{s}
			\end{equation}
			where $C_0=3 \frac{\max (\sqrt{d}, 3 \sqrt{2d}(s-1))^{2 s}}{s!}$. Then we have the nested sequence
			\begin{equation}\label{eq:sdp_relax}
			\V^{SDP}_{\sigma,rest}\subset \V^{AFF}_0\subset \V^{SDP}_0 \text{ for } \sigma:=\min(\sigma_{slow},\sigma_{fast}).
			\end{equation}
			If \Cref{ass:smooth_objective} also holds, then we have the bound
			\begin{equation}\label{eq:nested_optima_aposteriori}
			\Lcal(\bar{\b f}^{SDP}_{0,rest})\le \Lcal(\bar{\b f}^{AFF}_0)\le \overline{\Lcal}^{SDP}_{\sigma,rest},
			\end{equation}
			where $\overline{\Lcal}^{SDP}_{\sigma,rest}$ may be infinite, but is finite provided $\bar{\b f}^{AFF}_0 +\sigma \b g_1 \in \V^{SDP}_{0,rest}$. If $\Lcal(\cdot)$ is also $\beta$-Lipschitz on $\BB_K(\b 0, \bar{\b f}^{SDP}_{0,rest}+ \sigma \|\b g_1\|_K)$, with $\b g_1$ as in \Cref{thm:well-posed_cons}-\ref{it_nonempty}, then we have the bound
			\begin{equation}\label{eq:nested_optima_apriori}
			\Lcal(\bar{\b f}^{SDP}_{0,rest})\le \Lcal(\bar{\b f}^{AFF}_0)\le \Lcal(\bar{\b f}^{SDP}_{0,rest})+\beta \|\b g_1\|_{\HK} \sigma.
			\end{equation}

	\end{theorem}
    \noindent\tb{Discussion of the results of \Cref{thm:well-posed_cons,thm:main_result}.} By \Cref{thm:well-posed_cons}-\ref{it_aff_subset_eps} and the definition of $\V^{SDP}_{\sigma,rest}$ in \eqref{eq:def_sdp_constraints_rest}, it is clear that if there exists some $A^*\in S^+(\hh_\phi)$ such that $\b c(x)^\top \bar{\b f}^{AFF}_{0}(x)+d(x)= \scalh{\phi(x)}{A^* \phi(x)}$ for all $x\in\Kcons$ then $\bar{\b f}^{AFF}_{0}\in \V^{SDP}_{\sigma,rest}$ for $M_f(\bar{\b f}^{AFF}_0)=\|\bar{\b f}^{AFF}_0\|_K$ and $M_A(\bar{\b f}^{AFF}_0,\hat{X})=\tr(A^*)<\infty$. In every previous article on kSoS such as \citet{rudi2020finding,vacher2021ot}, the proof of fast convergence is divided into two steps: first the existence of $A^*$ is shown under some regularity assumptions, then a similar bound to \Cref{thm:main_result} is proven. We focused here on the second part, proposing a general approach to ensuring convergence. The latter occurs indeed at a faster rate $\sigma_{fast}$ if the problem exhibits some form of smoothness. Nevertheless the scheme always converges, at worst at rate $\sigma_{slow}$.

	\begin{proof}
		
		Let us show that \eqref{eq:sdp_relax} holds. Here we set $M_f=M_f(\bar{\b f}^{AFF}_0)$ and $M_A=M_A(\bar{\b f}^{AFF}_0,\hat{X})$ to shorten the notation and we deal separately with the cases $\sigma_{slow}$ and $\sigma_{fast}$ since the assumptions on $\Kcons$ and $\hat X$ differ.
		 
		 \tb{Case of $\sigma_{fast}$:} Take $\b f\in\V^{SDP}_{\sigma_{fast},rest}$ and $i\in[I]$. By definition of the constraint set, see \eqref{eq:def_sdp_constraints}, the traces of the corresponding $A_i\in S^+(\Hphi)$ are bounded by $M_A$. The extra assumptions on $\Kcons$ allow to apply \Cref{lem:scattered_inequalities} to $\xi_i(\cdot)=\b c_i(\cdot)^\top \b f(\cdot)+ d_i(\cdot)-\sigma_{fast}$ and $\tau=0$. Hence, $\xi_i(x) \geq-\epsilon_{0,i}$ for all $x \in \Omega$ with $\epsilon_{0,i}:=C_0\left(|\xi_i|_{\Omega, s}+2^s D_{\Omega,s}^2 M_A \right) \cdot (h_{\hat{X},\Kcons})^s$. Now we just have to upper bound $|\xi_i|_{\Omega, s}$. Let $\bm \alpha \in \N^d,\, |\bm \alpha|= s$. Since $|\xi_i|_{\Omega, s} \leq|\b c_i(\cdot)^\top \b f (\cdot)|_{\Omega,s}+\left|d\right|_{\Omega,s}$, the focus goes on $\p^{\bm \alpha}[\b c_i(\cdot)^\top \b f(\cdot)](x)$. By the reproducing property for derivatives \citep[Lemma 1][with a direct extension to non-constant coefficients $\b c(x)$]{aubin2020hard_SDP}, and using the definition of $C_{K,s}$ in \eqref{eq:err_SDP_s}
		 \begin{equation}\label{eq:bounding_variations_of_derivatives_constraints}
		 |\p^{\bm \alpha}[\b c_i(\cdot)^\top \b f(\cdot)](x)|=|\langle \b f,  \p_x^{\bm \alpha}[K(\cdot,x)\b c_i(x)]\rangle_{\HK}|\le C_{K,s} M_f.
		 \end{equation}
		 Consequently, $\xi_i(x) \geq- C_0\left(C_{K,s} M_f+\max_{i\in[I]}|d_i|_{\Omega, s}+2^s D_{\Omega,s}^2 M_A\right) h^{s}_{\hat{X},\Kcons}= -\sigma_{fast}$ for all $x \in \Omega$. Hence, by definition of $\xi_i$, $\b c_i(x)^\top \b f(x)+ d_i(x)\ge 0$, and $\b f \in \V^{AFF}_{0,rest}$. 
		 
		 \tb{Case of $\sigma_{slow}$:} In this case we just use a Lipschitzianity argument without applying \Cref{lem:scattered_inequalities}. Take $\b f\in\V^{SDP}_{\sigma_{slow},rest}$ and $i\in[I]$. Define $\xi_i$ as above, so $\xi_i(x_m)\ge 0$ for all $m\in[M]$, again by definition \eqref{eq:def_sdp_constraints}. Then, for any $x\in\Kcons$, there exists $m\in[M]$ such that $\|x-x_m\|\le h_{\hat{X},\Kcons}$. Applying \eqref{eq:bounding_variations_of_derivatives_constraints} for $|\bm \alpha|=1$, and then using that $\b f\in\V^{SDP}_{\sigma_{slow},rest}$, we obtain that
		 \begin{equation*}
		 \b c_i(x)^\top \b f(x)+ d_i(x) \ge \b c_i(x_m)^\top \b f(x)+ d_i(x_m) - h_{\hat{X},\Kcons} C_{K,1} M_f - h_{\hat{X},\Kcons} \max_{i\in[I]}|d_i|_{\Kcons, 1}\ge \sigma_{slow} - \sigma_{slow}= 0,
		 \end{equation*}
		 so $\b f \in \V^{AFF}_{0,rest}$.\\
		 
		As $\sigma=\min(\sigma_{slow},\sigma_{fast})$, we conclude that $\V^{SDP}_{\sigma,rest}$ is either equal to $\V^{SDP}_{\sigma_{slow},rest}$ or to $\V^{SDP}_{\sigma_{fast},rest}$, whence, by the proof above, $\V^{SDP}_{\sigma,rest}\subset \V^{AFF}_0$. By \Cref{thm:well-posed_cons}-\ref{it_aff_subset_eps}, $\V^{AFF}_0\subset \V^{SDP}_0$, so \eqref{eq:sdp_relax} holds. Putting the two inclusions together, we obtain that
		 \begin{equation}\label{eq:sdp_relax_restricted}
		 \V^{SDP}_{\sigma,rest}\subset \V^{AFF}_{0,rest}\subset \V^{SDP}_{0,rest}.
		 \end{equation}
		 By \Cref{thm:well-posed_cons}-\ref{it_solutions_existence}, we have the existence of $\bar{\b f}^{AFF}_0$ and of $\bar{\b f}^{SDP}_{0,rest}$. By \Cref{thm:well-posed_cons}-\ref{it_aff_subset_eps}, we have $\bar{\b f}^{AFF}_0\in \V^{AFF}_{0,rest}$ so it is also optimal for $\V^{AFF}_{0,rest}\subset \V^{AFF}_{0}$. We then directly deduce \eqref{eq:nested_optima_aposteriori} from the nestedness \eqref{eq:sdp_relax_restricted} of the constraint sets.\\
		 
		 We will now show \eqref {eq:nested_optima_apriori}, first proving that $\bar{\b f}^{SDP}_{0,rest} +\sigma\b g_1\in\V^{AFF}_{0}$.\footnote{Note that $\bar{\b f}^{SDP}_{0,rest} +\sigma\b g_1$ does not belong to $\V^{AFF}_{0,rest}$ in general due to the restriction of the constraint set by $M_f$ and $M_A$.} By definition of $\b g_1$, defining $\xi_i(\cdot)$ as follows, we get for any $x\in\Kcons$
		 \begin{align}\label{eq:f+g_Vaff}
		 \b c_i(x)^\top (\bar{\b f}^{SDP}_{0,rest} +\sigma\b g_1)(x) + d_i(x)=:\xi_i(x)  + \sigma \b c_i(x)^\top\b g_1(x)\ge \xi_i(x)  + \sigma.
		 \end{align}
		 With similar computations to the two cases discussed above, either applying a Lipschitzianity argument or \Cref{lem:scattered_inequalities}, we derive that $\xi_i(x)\ge -\sigma$. Hence $\bar{\b f}^{SDP}_{0,rest} +\sigma\b g_1\in\V^{AFF}_{0}$ for which $\bar{\b f}^{AFF}_0$ is optimal. Since $\Lcal(\cdot)$ is now also assumed to be $\beta$-Lipschitz on $\BB_K(\b 0, \bar{\b f}^{SDP}_{0,rest}+ \sigma \|\b g_1\|_K)$, we derive the bounds
		 \begin{equation}\label{eq:nested_optima_apriori2}
		 \Lcal(\bar{\b f}^{SDP}_{0,rest})\le \Lcal(\bar{\b f}^{AFF}_0)\le \Lcal(\bar{\b f}^{SDP}_{0,rest}+\sigma\b g_1)\le \Lcal(\bar{\b f}^{SDP}_{0,rest})+\beta \|\b g_1\|_{\HK} \sigma.
		 \end{equation}
		 Note that in \eqref{eq:sdp_relax_restricted}, $\V^{SDP}_{\sigma,rest}$ may be empty as the restriction depends on $\bar{\b f}^{AFF}_0$ rather than $\bar{\b f}^{AFF}_{\sigma}$ (the bounds on the functions may be too tight). Consequently $\overline{\Lcal}^{SDP}_{\sigma,rest}$ may be infinite. If $\bar{\b f}^{AFF}_0 +\sigma \b g_1 \in \V^{SDP}_{0,rest}$, then, proceeding as in \eqref{eq:f+g_Vaff}, we can show that $\bar{\b f}^{AFF}_0 +\sigma \b g_1 \in \V^{AFF}_{\sigma,rest}\subset \V^{SDP}_{\sigma,rest}$, which implies that $\overline{\Lcal}^{SDP}_{\sigma,rest}$ is finite as the constraint set is non-empty and since  \Cref{ass:smooth_objective} was assumed.
	\end{proof}

	\section{Numerical application}\label{sec:numerics}
	
	We illustrate our methodology on the problem of learning vector fields under constraints, since the two other examples mentioned in \Cref{sec:examples} were already showcased in earlier works \citep{rudi2020finding,vacher2021ot}. Inspired by the example in \citet{ahmadi23learning}, we consider a competitive Lotka-Volterra model in two dimensions over the set $\X:=[0,1]^2$,
	\begin{equation}\label{eq:example_GenLotka}
	\dot{\b x}(t)=\bar{\b f}(\b x(t)), \text{ where } \b x(t) \in \mathbb{R}^{2}\text{ and } \bar{\b f}(\b x)=\begin{pmatrix}
	x_{1}\left(1-x_{1} -0.2\cdot x_{2}\right)\\
	x_{2}\left(1-x_{2} -0.4\cdot x_{1}\right)
	\end{pmatrix}
	\end{equation}
	The structure of the model incites to introduce the following interpolation and invariance constraints, for $\Kcons:=\partial \X$,
	\begin{align}
		\b f(\b 0)&=\b 0; \label{eq:GenLotka_interp} \tag{Interp}\\
		\scalidx{\overrightarrow{\b n}(\b x)}{\b f(\b x)}{\R^2} &\le 0, \, \forall \, \b x\in \Kcons; \label{eq:GenLotka_inv} \tag{Inv}
	\end{align}
	where $\overrightarrow{\b n}(\b x)$ is the outer-pointing normal to $\Kcons=\partial ([0,1]^2)$ at point $\b x$ (we consider two normals on corners). As discussed above, we discretize the inequality constraint at $M$ points $\{\tilde{\b x}_m\}_{m\in[M]}\subset \Kcons$ (with possible repetition), defining the sampled and kSoS invariance constraints (for some $A\in S^+(\Hphi)$):
	\begin{align}
	\scalidx{\overrightarrow{\b n}(\tilde{\b x}_m)}{\b f(\tilde{\b x}_m)}{\R^2} &\le 0, \, \forall \, m\in[M]; \label{eq:GenLotka_sampled-inv} \tag{Sampled-Inv}\\
	-\scalidx{\overrightarrow{\b n}(\tilde{\b x}_m)}{\b f(\tilde{\b x}_m)}{\R^2} &= \scalidx{\phi(\tilde{\b x}_{m})}{A \phi(\tilde{\b x}_{m})}{\Hphi}, \, \forall \, m\in[M]; \label{eq:GenLotka_kSoSinv} \tag{kSoS-Inv}
	\end{align}
	Suppose we are given a few noisy snapshots forming a training dataset
	\begin{equation*}
		\mathcal{D}:=\left\{\left(\b x_n, \b y_n\right)\right\}_{n\in[N]}=\left\{\left(\b x_n, \bar{\b f}(\b x_n)\right)+10^{-2}\cdot(\bm \varepsilon_{n, 1},\bm \varepsilon_{n, 2})\right\}_{n\in[N]}
	\end{equation*}
		where the noises $\bm \varepsilon_{n, j}$ are independently sampled from the standard normal distribution. We then solve the kernel ridge regression (KRR) problem
	\begin{mini}|s|
		{\substack{\b f \in \HK,\\ A\in S^+(\hh_\phi)\\+ \text{constraints}}}{\sum_{n\in[N]} \|\b y_n-\b f(\b x_n)\|^2_{\R^2}+\la_K \|\b f\|_K^2 + \la_\phi \tr(A) }{\label{opt-cons_KRR}}{}%
	\end{mini}
	either without constraints ($\la_\phi=0$), or with sampled constraints ($\la_\phi=0$, $\b f$ satisfying \eqref{eq:GenLotka_interp} and \eqref{eq:GenLotka_sampled-inv}), or with kSoS-sampled constraints ($\la_\phi>0$, $\b f$ satisfying \eqref{eq:GenLotka_interp} and \eqref{eq:GenLotka_kSoSinv}).
	
	Denoting the Gaussian kernel of bandwidth $\sigma_K>0$ by $k_{\sigma_K}(\b x, \b y)=\exp\left(\frac{\text{-}\|\b x - \b y\|^2_{\R^2}}{2 \sigma_K^2}\right)$, we choose $K$ to be a diagonal Gaussian Kernel, $K(\b x, \b y)=k_{\sigma_K}(\b x, \b y)\Id_2$, and $k_{\phi}(\b x, \b y)=k_{\sigma_\phi}(\b x, \b y)$ for some $\sigma_K, \sigma_\phi>0$. We thus have up to four hyperparameters to select. We start with $\sigma_K$ that is optimized by gridsearch over $[0.1,2]$ for $\la_K$ deterministically obtained by generalized cross-validation \citep{Golub1979}. We then use the heuristic $\sigma_\phi=\sigma_K/2$, motivated by the sum-of-squares formulation of \eqref{eq:GenLotka_kSoSinv}, and fix $\la_\phi$ by gridsearch over the logarithmic interval $[10^{-6}, 10^{-1}]$. In the experiment we set $N=5$ and draw $\{\b x_n\}_{n\in[N]}$ uniformly at random in $\X$. The $M$ sample points are taken on a regular grid on $\Kcons$, consequently there are $M/4$ constraints per side of $\X$. The optimal hyperparameters resulted in $\sigma_K=1$ and $\la_\phi=10^{-3}$. In \Cref{sec:theory}, we studied constraints over $\tr(A)$, however these require to have a good estimate beforehand of the bound $M_A$, facing the risk of the problem having no solution for $M_A$ too small. The penalization hyperparameter $\la_\phi$ is instead much easier to tune, and \eqref{opt-cons_KRR} always has a solution. \rv{The implementation on a laptop was done on Matlab through the code generator CVXGEN, using SDPT3 as solver. The CPU times are reported in \Cref{table:cpu_times}. Off-the-shelf solvers notoriously struggle on SDP problems, the computation time of the kernel matrices being negligible in comparison. Nevertheless the specific structure of kSoS constraints is amenable for more specific algorithms such as the Damped Newton scheme used in \citep[Section 6]{rudi2020finding}. One will still pay a price for dealing with SDP constraints, but not as dire.}
	
\begin{figure}[ht!]
		\centering	
		\begin{minipage}[c]{.4\linewidth}
			\begin{center}			
				\subfloat[][True vector
				field]{\resizebox{\linewidth}{!}{\label{VectorField_GenLotka_true}\includegraphics{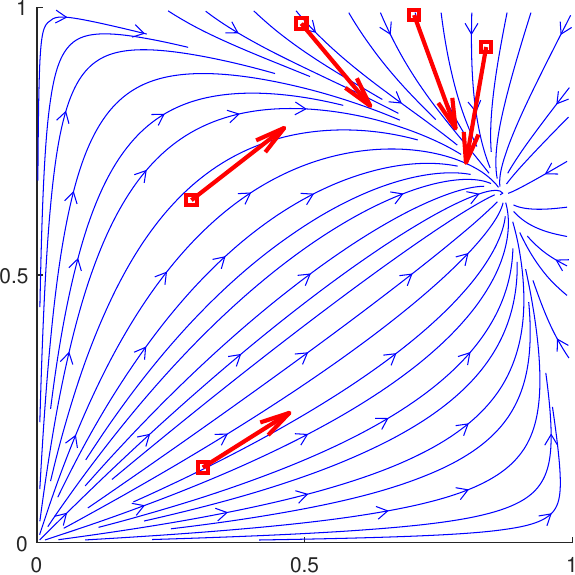}}}
			\end{center}				
		\end{minipage}
		\begin{minipage}[c]{.4\linewidth}
			\begin{center}			
				\subfloat[][Unconstrained KRR]{\resizebox{\linewidth}{!}{\label{VectorField_GenLotka_NoCons}\includegraphics{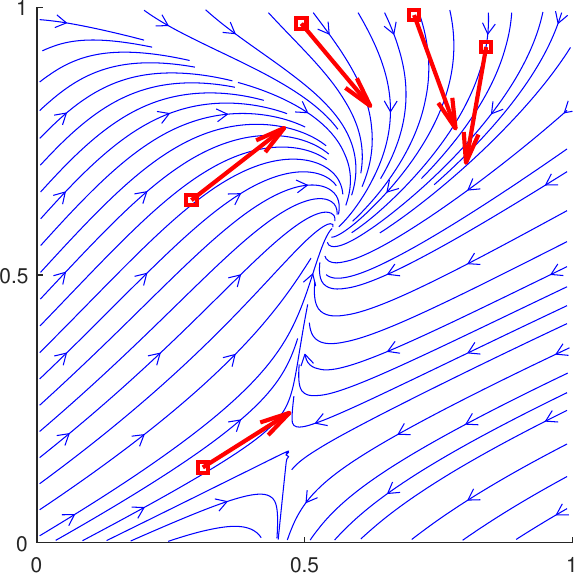}}}
			\end{center}				
		\end{minipage}
		
		\begin{minipage}[c]{.4\linewidth}
			\begin{center}			
				\subfloat[][Interp $\cap$ Sampled-Inv]{\resizebox{\linewidth}{!}{\label{VectorField_GenLotka_0Inv}\includegraphics{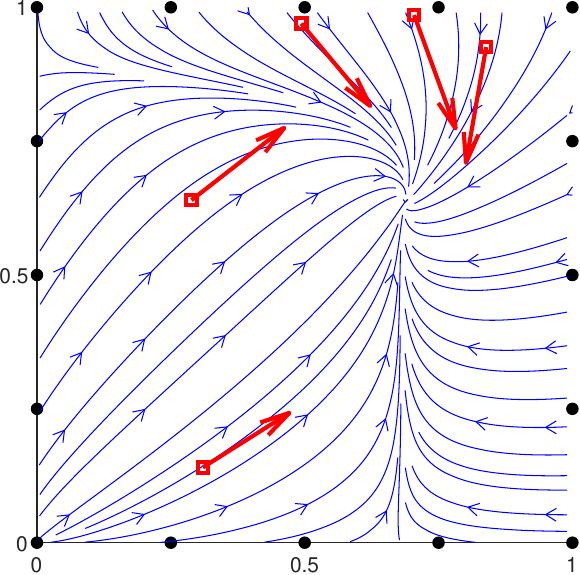}}}
			\end{center}				
		\end{minipage}
		\begin{minipage}[c]{.4\linewidth}
			\begin{center}			
				\subfloat[][Interp $\cap$ kSoS-Inv]{\resizebox{\linewidth}{!}{\label{VectorField_GenLotka_0kSoSInv}\includegraphics{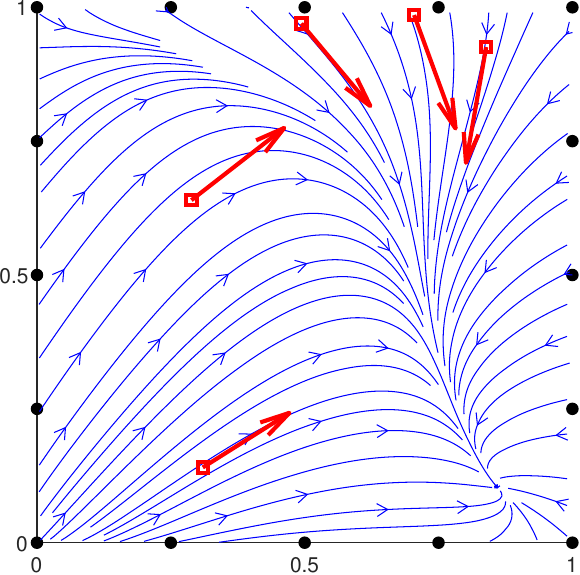}}}
			\end{center}				
		\end{minipage}
		\caption{Streamplots of the vector field in \eqref{eq:example_GenLotka} (\Cref{VectorField_GenLotka_true}) and of kernel solutions with different constraints (\Cref{VectorField_GenLotka_NoCons,VectorField_GenLotka_0Inv,VectorField_GenLotka_0Inv}). Red arrows: $N=5$ noisy observations. Black points: $M=20$ samples on the constraint set $\Kcons$.}
		\label{fig:VectorField_GenLotka}
	\end{figure}

\begin{figure}[ht!]
		\centering	
		\begin{minipage}[c]{.4\linewidth}
			\begin{center}			
				\subfloat[][Angular error]{\resizebox{\linewidth}{!}{\label{CosErr_GenLotka_true}\includegraphics{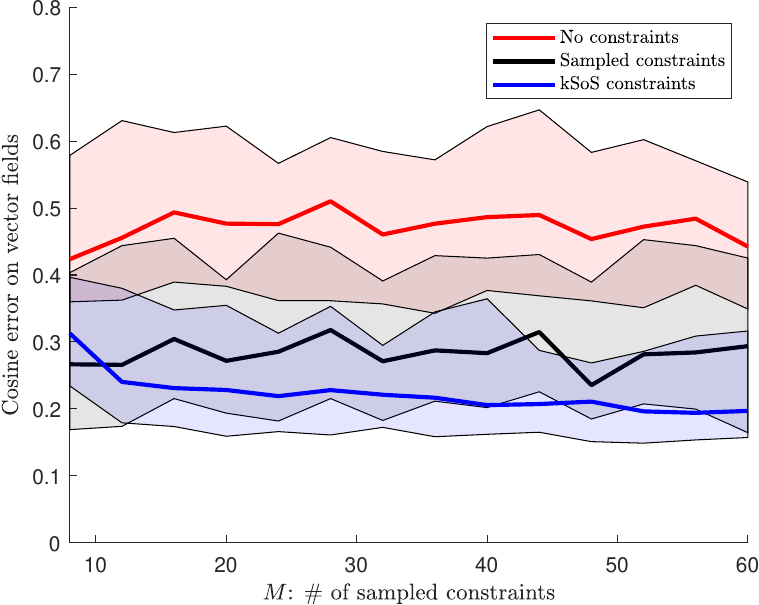}}}
			\end{center}				
		\end{minipage}
		\begin{minipage}[c]{.4\linewidth}
			\begin{center}			
				\subfloat[][$L^2$  error]{\resizebox{\linewidth}{!}{\label{L2Err_GenLotka_NoCons}\includegraphics{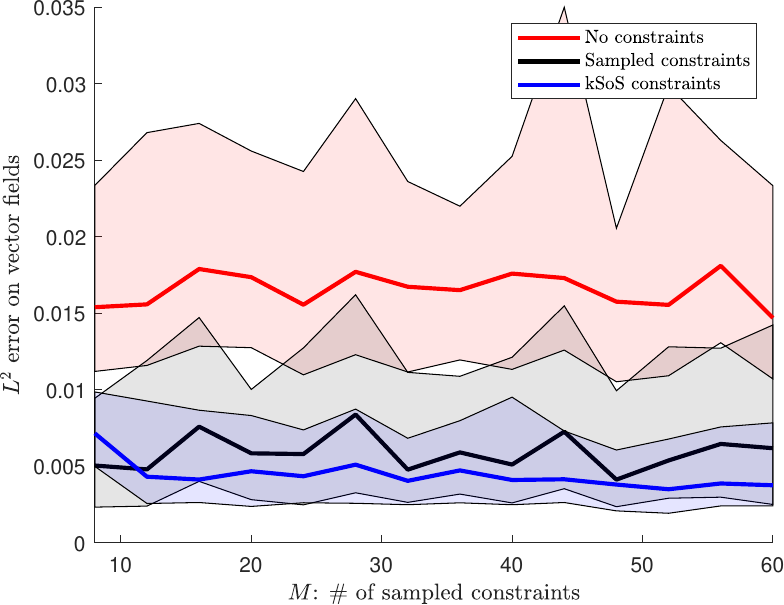}}}
			\end{center}				
		\end{minipage}
		
		\begin{minipage}[c]{.4\linewidth}
			\begin{center}			
				\subfloat[][$L^\infty$ violation on $\Kcons$]{\resizebox{\linewidth}{!}{\label{LinfErr_GenLotka_0Inv}\includegraphics{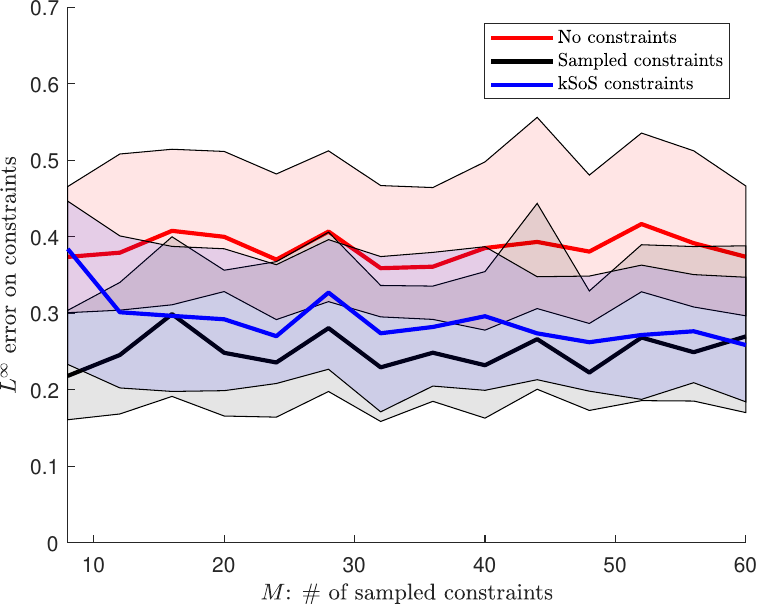}}}
			\end{center}				
		\end{minipage}
		\begin{minipage}[c]{.4\linewidth}
			\begin{center}			
				\subfloat[][$L^1$ violation on $\Kcons$]{\resizebox{\linewidth}{!}{\label{L1Err_GenLotka_0kSoSInv}\includegraphics{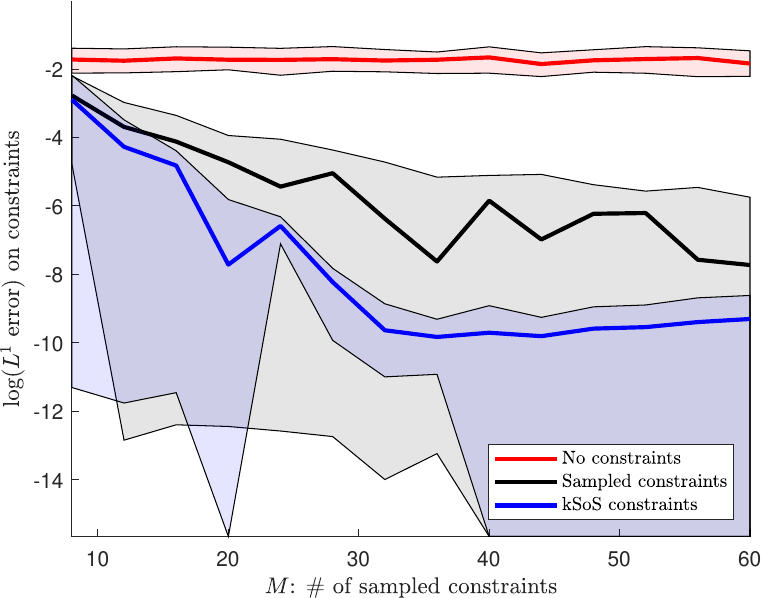}}}
			\end{center}				
		\end{minipage}
		\caption{Performance criteria for the reconstruction of the vector field in \eqref{eq:example_GenLotka} as a function of the number $M$ of inequality constraints imposed}
		\label{fig:ErrPlot_GenLotka}
	\end{figure}

	We report the performance of our reconstruction on \Cref{fig:ErrPlot_GenLotka} based on the following four criteria: i) an angular measure of error $\frac{1}{|\X|}\int_\X\arccos\left(\frac{\scalidx{\hat{\b f}(x)}{\bar{\b f}(x)}{\R^P}}{\|\hat{\b f}(x)\|\cdot \|\bar{\b f}(x)\|}\right)\d x$ which reflects that the phase portraits as in \Cref{fig:VectorField_GenLotka} are similar since it disregards the magnitude of the fields; ii) the $L^2$ reconstruction error $\frac{1}{|\X|}\int_X\|\hat{\b f}(x)- \bar{\b f}(x)\|^2\d x$; iii) the $L^\infty$ measure of violation of constraints $\text{-}\min_{x\in\Kcons,i\in\Ical}(0,\b c_i(x)^\top\hat{\b f}(x)+d_i(x))$; iv) the $L^1$ measure of violation of constraints $\frac{\text{-}1}{|\Kcons|}\int_\Kcons\min_{i\in\Ical}(0,\b c_i(x)^\top\hat{\b f}(x)+d_i(x))\d x$. All the integrals are approximated by a uniform grid with a 100 points per dimension. The first quartile, the median and the third quartile of the values computed are obtained by repeating a 100 times each experiment when varying the number $M$ of sampled constraints.
	
 	On \Cref{VectorField_GenLotka_NoCons}, we observe that the solution of unconstrained KRR fails to satisfy the invariance assumption \eqref{eq:GenLotka_inv}. On the contrary, both the constrained solutions satisfy it. Moreover the constraints \eqref{eq:GenLotka_kSoSinv} provides a streamplot (\Cref{VectorField_GenLotka_0kSoSInv}) looking more alike the true $\bar{\b f}$ (\Cref{VectorField_GenLotka_true}) than the result obtained through sampled constraints \eqref{eq:GenLotka_sampled-inv} (\Cref{VectorField_GenLotka_0Inv}). This is further ascertained numerically on \Cref{fig:ErrPlot_GenLotka} where the kSoS constraint has a better reconstruction error in both angular and $L^2$ measures (\Cref{CosErr_GenLotka_true,L2Err_GenLotka_NoCons}). \Cref{LinfErr_GenLotka_0Inv} shows that the constraint \Cref{eq:GenLotka_inv} is nevertheless violated in most experiments with either \eqref{eq:GenLotka_kSoSinv} or \eqref{eq:GenLotka_sampled-inv}. However \Cref{L1Err_GenLotka_0kSoSInv} illustrates that the $L^1$ measure of violation decreases faster for the kSoS sampling w.r.t.\ the sampled constraints. The kSoS constraint indeed leverages the smoothness of the functions. If the constraints were satisfied in an experiment, the error was then set to machine precision to be able to take the logarithm, since the corresponding violation was null.

    \begin{table}[!h]
    \small
\centering
\rv{\begin{tabular}{l|cccccc}
M (\# samples of constraints) & 8 & 48 & 88 & 128 & 168 & 208 \\
\midrule
CPU time (s) for AFF & 1.2 $\pm$ 6e-1 & 2.2 $\pm$ 9e-1 & 1.7 $\pm$ 8e-1 & 3.9 $\pm$ 1 & 5.3 $\pm$ 6e-1 & 6.2 $\pm$ 1 \\
CPU time (s) for SDP & 1.1 $\pm$ 2e-1 & 2.2 $\pm$ 6e-1 & 5.6 $\pm$ 5e-1 & 18.3 $\pm$ 2 & 74.2 $\pm$ 8 & 468.8 $\pm$ 31 \\
\end{tabular}
\caption{CPU times in seconds, comibining the code generation and solver for solving the quadratic problem either with the affine (discretized) constraints (AFF) or with the kSoS constraints (SDP). SDP constraints are notoriously harder to solve than QPs with off-the-shelf solvers.}\label{table:cpu_times}}

\end{table}

	\section{Technical results on closedness of sets and on scattering inequalities}\label{sec:closed-sets}

	\begin{lemma}[Closed constraint sets]\label{lem:closed_sdp_constraints} Fix $\epsArr\in\R^I$, $M_f\ge 0$ and $M_A\ge 0$. Then $\V^{AFF}_{0}$ is a weakly closed convex subset of $\HK$, and, if \Cref{ass:constraint_set} also holds, so are $\V^{SDP}_{\epsilon}$ and $\V^{SDP}_{\epsilon,rest}$.
	\end{lemma}
	\begin{proof}
		The sets $\V^{AFF}_{0}$ and $\V^{SDP}_{\epsilon}$ are clearly convex, so is $\V^{SDP}_{\epsilon,rest}$ using the triangular inequality with the nuclear/trace norm on $S^+(\Hphi)$. Using Mazur's lemma \citep[see e.g.\ Proposition 1][p34]{beauzamy1985introduction}, it then suffices to show that the sets are strongly closed to conclude that they are weakly closed in $\HK$. For $\V^{SDP}_{\epsilon}$, the proof is quite straightforward provided \Cref{ass:constraint_set} is satisfied. Indeed, fix $(\b f_n)_{n\in \N}\in (\V^{SDP}_{\epsilon})^\N$ converging to $\b f\in \HK$. As $\HK$ is a RKHS, for any $m\in[M]$, $\b y_{m,n}:=\b C (x_m) \b f_n (x_m)+ \b d(x_m)-\epsArr$ converges to $\b y_{m}:=\b C (x_m) \b f (x_m)+ \b d(x_m)-\epsArr$. Since $\b f_n\in \V^{SDP}_{\epsilon}$, we have that $\b y_{m,n}\ge 0$, so $\b y_{m}\ge 0$ and \Cref{lem:kSoS_interpolation} gives the $(A_i)_{i\in[I]}$ corresponding to $\b f$ and $\b f\in\V^{SDP}_{\epsilon}$. Hence $\V^{SDP}_{\epsilon}$ is strongly closed in $\HK$. The same nonegativity argument paired with the fact that $\HK$ is a RKHS gives the strong closedness of $\V^{AFF}_{0}$. 
		
		For $\V^{SDP}_{\epsilon,rest}$ instead, we have to bound $\tr(A_i)$ by $M_A$. For this, we will have to explicitly relate $\V^{SDP}_{\epsilon,rest}$ to a compact set in finite dimensions. Similarly to the proof of \citet[Lemma 2]{rudi2020finding} and to that of \Cref{lem:kSoS_interpolation} above, let $\widehat{S}:\Hphi \rightarrow \RR^M$ be the linear operator acting as follows
		$$
		\widehat{S} g=\left(\left\langle\phi(x_{1}), g\right\rangle_{\Hphi}, \ldots,\left\langle\phi(x_{M}), g\right\rangle_{\Hphi}\right) \in \RR^M, \quad \forall g \in \Hphi
		$$
		Consider the adjoint of $\widehat{S}$, $\widehat{S}^{*}: \RR^M \rightarrow \Hphi$, defined consequently as $\widehat{S}^{*} \beta=\sum_{m\in[M]} \beta_m \phi\left(x_m\right)$ for $\bm \beta \in \RR^M$. Let $\b K_{\phi,\hat{X}}=[k_\phi(x_i,x_j)]_{i,j\in[M]}=\b  R^\top  \b R$, where $\b R:=[\Phi_1,\dots, \Phi_M]$ is the Cholesky decomposition of the symmetric $\b K_{\phi,\hat{X}}\in\RR^{M\times M}$ and $\Phi_m=\b R \b e_m$ are its columns. Note, in particular, that $\b K_{\phi,\hat{X}}=\widehat{S} \widehat{S}^{*}$ and that $\widehat{S}^{*} e_{m}=\phi\left(x_m\right)$, where $e_{m}$ is the $m$-th element of the canonical basis of $\mathbb{R}^{M}$.
		
		As $\b K_{\phi,\hat{X}}$ is invertible by \Cref{ass:constraint_set}, we define the operator $V=\b R^{-\top} \widehat{S}$ and its adjoint $V^{*}=\widehat{S}^{*} \b R^{-1}$. By using the definition of $V$, that $\b K_{\phi,\hat{X}}=\b R^{\top} \b R$, and that $\b K_{\phi,\hat{X}}=\widehat{S} \widehat{S}^{*}$, we derive two facts.
		On the one hand,
		$$
		V V^{*}=\b R^{-\top} \widehat{S} \widehat{S}^{*} \b R^{-1}=\b R^{-\top} \b K_{\phi,\hat{X}} \b R^{-1}=\b R^{-\top} \b R^{\top} \b R \b R^{-1}=\text{Id}_{\RR^M}
		$$
		On the other hand, $P:=V^{*} V:\HK\rightarrow\HK$ is a projection operator, i.e., $P^{2}=P, P$ is positive definite and its range is $\Sp\left\{\phi\left(x_m\right) \mid m \in[M]\right\}$, implying $P \phi\left(x_m\right)=\phi\left(x_m\right)$ for all $m \in[M] .$ Indeed, using the equation above, $P^{2}=V^{*} V V^{*} V = V^{*}\left(V V^{*}\right) V=V^{*} V=P$, and the positive semi-definiteness of $P$ is given by construction since it is the product of an operator and its adjoint. Moreover, the range of $P$ is the same as that of $V^{*}$ which in turn is the same as that of $S^{*}$, since $\b R$ is invertible. 
		
		Let $C_{\hat{X}}:\HK \rightarrow \RR^M$ be the linear operator defined as follows
		\begin{equation*}
		C_{\hat{X}} \b f=\left(\b C(x_1) \b f(x_1),\dots, \b C(x_M) \b f(x_M)\right) \in \RR^{I\times M}, \quad \forall\, \b f \in \HK.
		\end{equation*}
		Consider the subspace $ \E:=C_{\hat{X}}(\HK)\subset \RR^{I\times M}$. From now on, we restrict the output space, i.e.\ we focus on the surjective $C_{\hat{X}}:\HK \rightarrow \E$. We may thus introduce its injective right inverse $C_{\hat{X}}^+:\E \rightarrow \HK$ defined as
		\begin{equation*}
		C_{\hat{X}}^+ \b V= \argmin{\substack{\b f\in\HK,\\ C_{\hat{X}} \b f= \b V}}\| \b f \|_{\HK}, \quad \forall\, \b V \in \E.
		\end{equation*}
		The continuous operator $C_{\hat{X}}^+$ is related to a subspace of $\HK$ since $\HK=\Ker(C_{\hat{X}})\oplus C_{\hat{X}}^+(\E)$. As $C_{\hat{X}}^+$ is linear, $C_{\hat{X}}^+(\E)$ is finite dimensional. Fix the compact set
		\begin{align*}
		\mathcal{R}_{K,\phi}=\{\b V \in \E\, | \, 	C_{\hat{X}}^+ \b V\in \BB_K(0,M_f),\, \exists &(\b B_i)_{i\in[I]}\in S^+(R^M)^I,\, \tr(\b B_i)\le M_A,\,\\ &\b V \b e_m + \b d(x_m)=\epsArr+[\Phi_m^\top \b B_i \Phi_m]_{i\in[I]},\, \forall \, m\in [M] \}.
		\end{align*}
		For any $\b f\in \V^{SDP}_{\epsilon,rest}$ and $\b g\in\Ker(C_{\hat{X}})$, as $C_{\hat{X}} (\b f+\b g)=C_{\hat{X}}\b f$, the values appearing in the definition of $\V^{SDP}_{\epsilon,rest}$ do not change, so $\b f+\b g\in \V^{SDP}_{\epsilon,rest}$. Hence, setting $\G_{rest}:=C_{\hat{X}}^+(\E)\cap \V^{SDP}_{\epsilon,rest}$, we have that $\V^{SDP}_{\epsilon,rest}=\G_{rest} \oplus \Ker(C_{\hat{X}})$.
		
		We now show that $\G_{rest}=C_{\hat{X}}^+(\mathcal{R}_{K,\phi})$. Let $\b V\in \mathcal{R}_{K,\phi}$, and consider the corresponding $(\b B_i)_{i\in[I]}\in S^+(R^M)^I$. Let $A_i:=V^* \b B_i V\in S^+(\Hphi)$, then, by the cyclicity of the trace and since $V V^*=\Id$, $\tr(A_i)=\tr(\b B_i V V^*)=\tr(\b B_i)\le M_A$. So $C_{\hat{X}}^+(\b V)\in C_{\hat{X}}^+(\E)\cap \V^{SDP}_{\epsilon,rest}=\G_{rest}$. Conversely, let $\b f \in \G_{rest}$ and consider the corresponding $(A_i)_{i\in[I]}\in S^+(\Hphi)^I$. Then $\b B_i:=V \b B_i V^*\in S^+(R^M)$ satisfies $\tr(\b B_i)=\tr(\b B_i V^* V)\le \|V^* V\|_{op} \tr(A_i)\le M_A$ as $V^* V$ is a projector and, for any operators $B, C\in S^+(\Hphi)$, $\tr(BC)\le \|B\|_{op} \tr(C)$ \citep[e.g.][Lemma 1]{sheng1986trace}. So $C_{\hat{X}}\b f\in \mathcal{R}_{K,\phi}$, and, as $C_{\hat{X}}^+ $ is injective, $\b f\in C_{\hat{X}}^+(\mathcal{R}_{K,\phi})$. We thus conclude that $\G_{rest}=C_{\hat{X}}^+(\mathcal{R}_{K,\phi})$, thus $\G_{rest}$ is a compact set for the strong topology of $\HK$. So $\V^{SDP}_{\epsilon,rest}=\G_{rest} \oplus \Ker(C_{\hat{X}})$ is strongly closed in $\HK$, which concludes the proof.
	\end{proof}
	
	\citet{rudi2020finding} have focused mostly on multiplicative algebra\footnote{A vector space $\Hk$ of real-valued functions is a multiplicative algebra if there exists $c\in\R_+$ such that, for all $f,g \in \Hk$, the pointwise product $f\cdot g$ belongs to $\Hk$, and $\|f\cdot g\|\le c \|f\| \|g\|$}, among which Sobolev spaces. Here instead we tackle more general RKHSs. 
	
	\begin{lemma}[{{Modified \citet[][Theorem 4, p14]{rudi2020finding}}}]\label{lem:scattered_inequalities}
		Let $\X=\R^d$, take $s_\phi\ge 1$ and $s$ such that $1\le s \le s_\phi$. Consider a bounded set $S\subset \X$ and some $r>0$, set $\Omega:= \cup_{x\in S}\dBB(x,r)$.\footnote{Scattering inequalities such as \citet[][Theorem 4, p14]{rudi2020finding} unfortunately require to prescribe the form of $\Kcons$. If $\Kcons$ is any compact, we cannot get away with just a dilation. Indeed if $\hat{X}\subset \Kcons$ with $\Kcons+\dBB(0,r) \subset \cup_{m\in[M]} \BB_{\R^d}(x_m,\delta_m)$ then we face a problem, because $r\le\|\bm \delta\|_\infty$ but $r\le\|\bm \delta\|_\infty\le r \min(1,\frac{1}{18(s-1)^2}) $ cannot hold.} Take $\hat{X}=\{x_m\}_{m\in[M]}\subset \Omega$ with $h_{\hat{X}, \Omega}\le r \min(1,\frac{1}{18(s-1)^2}) $. Let $g\in C^s(\X,\R)$ and assume there exists $A\in S^+(\Hphi)$ and $\tau \ge 0$ such that
		\begin{equation}\label{eq:approx_by_sos}
		|g(x_m)-\langle \phi(x_m),  A \phi(x_m)\rangle_{\Hphi}|\le \tau, \quad \forall\, m\in[M],
		\end{equation}
		then the following statement holds
		\begin{equation}\label{eq:violation_amount}
		g(x) \geq-(\epsilon+2 \tau)\quad \forall x \in \overline{\Omega}, \quad \text { where } \epsilon=C (h_{\hat{X}, \Omega})^s
		\end{equation}
		and $C=C_{0}\left(|g|_{\Omega, s}+2^s D_{\Omega,s}^2 \tr(A)\right)$ with $C_{0}=3 \frac{\max (\sqrt{d}, 3 \sqrt{2d}(s-1))^{2 s}}{s!}$.
	\end{lemma}
	\begin{proof} 
		The proof is mostly similar to that of \citet[][Theorem 4, p14]{rudi2020finding}, also relying heavily on the technical \citet[][Theorem 13, p39]{rudi2020finding}, but discarding \citet[][Lemma 9, p41]{rudi2020finding} which requires multiplicative algebra. Define the functions $u, r_{A}: \Omega \rightarrow \R$ as follows
		$$
		r_{A}(x)=\langle\phi(x), A \phi(x)\rangle, \quad u(x)=g(x)-r_{A}(x), \quad \forall x \in \Omega
		$$
		By assumption, $\left|u\left(x_m\right)\right| \leq \tau$ for any $m \in[M]$. Since $h_{\bar{X}, \Omega} \le \|\bm \delta\|_\infty$, $\Omega:=\cup_{x\in\Kcons}\dBB(x,r)$, and we have assumed  $\|\bm \delta\|_\infty\leq r / \max \left(1,18(m-1)^{2}\right)$, we can apply apply \citet[][Theorem 13, p39]{rudi2020finding} (related to classical results on functions with scattered zeros \citep[][Chapter 11]{wendland2005scattered}) to bound  $ \sup_{x \in \Omega}|u(x)|$, whence
		$$
		\sup_{x \in \Omega}|u(x)| \leq 2 \tau+\epsilon, \quad \epsilon=c R_{s}(u) h_{\bar{X}, \Omega}^{s}
		$$
		where $c=3 \max \left(1,18(s-1)^{2}\right)^{s}$ and $R_{s}(v)=\sum_{|\bm \alpha|=s} \frac{1}{\bm \alpha!} \sup_{z \in \Omega}\left|\p^{\bm \alpha} v(x)\right|$ for any $v \in$ $C^{s}(\Omega)$. Since $r_{A}(x)=\langle\phi(x), A \phi(x)\rangle \geq 0$ for any $x \in \Omega$ as $A \in S^{+}(\Hphi)$, we have that 
		\begin{equation}\label{eq:lowerbound_eps_proof_scattering}
		g(x) \geq g(x)-r_{A}(x)=u(x) \geq-|u(x)| \geq-(2 \tau+\varepsilon), \quad \forall x \in \Omega	
		\end{equation}
		The last step is bounding $R_{s}(u)$. Recall the definition of $|u|_{\Omega,s}$ from \eqref{eq:deriv_bound}. Since $A$ is trace class, it admits a singular value decomposition $A=\sum_{i \in \N} \sigma_{i} u_{i} \otimes v_{i} .$ Here, $\left(\sigma_m\right)_{j \in \N}$ is a non-increasing sequence of non-negative eigenvalues converging to zero, and $\left(u_m\right)_{j \in \N}$ and $\left(v_m\right)_{j \in \N}$ are two orthonormal families of corresponding eigenvectors, (a family $\left(e_m\right)$ is said to be orthonormal if for $i, j \in \N,\left\langle e_{i}, e_m\right\rangle=1$ if $i=j$ and $\left\langle e_{i}, e_m\right\rangle=0$ otherwise). Note that we can write $r_{A}$ using this decomposition as, by the reproducing property, $r_{A}(x)=\sum_{i \in \N} \sigma_{i} u_{i}(x). v_{i}(x)$ with $u_{i},v_{i}\in\HK\subset \C^{s_\phi}(\X,\R)$.%
		
		Since $k_\phi(\cdot,\cdot) \in \C^{s_\phi,s_\phi}\left(\X\times \X,\R\right)$, we have that $u_m, v_m\in \C^{s_\phi}\left(\X,\R\right)$, so we can take derivatives of $r_A$. By Leibniz's formula extended to partial derivatives of products, also using the reproducing property for derivatives \citep[Lemma 1]{aubin2020hard_SDP} and the Cauchy-Schwarz inequality, for any $\bm \alpha\in \N^d$ with $|\bm\alpha|\le s$, denoting by $\mathcal{P}([\bm \alpha])$ the set of parts of $\iv{0}{\alpha_1}\times \dots \iv{0}{\alpha_d}$
		\begin{align*}
		\p^{\bm \alpha} r_{A}(x)=\sum_{i \in \N} \sigma_{i} \p^{\bm \alpha} (u_{i}\cdot v_{i})(x)&=\sum_{i \in \N} \sigma_{i} \sum_{S \in \mathcal{P}([\bm \alpha])} \frac{\p^{|S|}u_{i}}{\prod_{j\in S} \p x_j}(x) \frac{\p^{|\bm \alpha|-|S|}v_{i}}{\prod_{j\notin S} \p x_j}(x)\\
		& \le \sum_{i \in \N} \sigma_{i} \sum_{S \in \mathcal{P}([\bm \alpha])} \langle u_{i}, \frac{\p^{|S|}k_\phi(x,\cdot)}{\prod_{j\in S}\p x_j} \rangle_{\Hphi}\langle v_{i}, \frac{\p^{|\bm \alpha|-|S|}k_\phi(x,\cdot)}{\prod_{j\notin S}\p x_j} \rangle_{\Hphi}\\
		&\le \sum_{i \in \N} \sigma_{i} 2^s D_{\Omega,s}^2 \|u_i\|_{\Hphi} \|v_i\|_{\Hphi}.
		\end{align*}
		Since $\|u_i\|_{\Hphi}=\|v_i\|_{\Hphi}=1$, we obtain that $|r_A|_{\Omega,s}\le \tr(A) 2^s D_{\Omega,s}^2$. To conclude, note that, by the multinomial theorem,
		$$
		R_{s}(u)=\sum_{|\bm \alpha |=s} \frac{1}{\bm \alpha !} \sup_{x \in \Omega}\left|\p^{\bm \alpha} u(x)\right| \leq \sum_{|\bm \alpha |=s} \frac{1}{\bm \alpha  !}|u|_{\Omega,s}=\frac{d^{s}}{s !}|u|_{\Omega, s}.
		$$
		Since $|u|_{\Omega, s} \leq|g|_{\Omega,s}+\left|r_{A}\right|_{\Omega,s}$, combining all the previous bounds, we obtain
		$$
		\varepsilon \leq C_{0}\left(|g|_{\Omega, s}+ 2^s D_{\Omega,s}^2 \tr(A)\right) h_{\hat{X}, \Omega}^{s}, \quad \text{with} \quad C_{0}=3 \frac{d^{s} \max \left(1,18(s-1)^{2}\right)^{s}}{s !}.
		$$
		The proof is concluded by bounding $\epsilon$ in \eqref{eq:lowerbound_eps_proof_scattering} with the inequality above. Since $g\in C^0(\X,\R)$, \eqref{eq:lowerbound_eps_proof_scattering} also holds for $x \in \overline{\Omega}$ by continuity.

	\end{proof}

	\section{Appendix: Selection theorem and lower semicontinuity of a set-valued map}\label{sec:set-valued_lsc}
	In this section we return to \eqref{opt-cons_gen} without the affine subspace constraint, and study the existence of a solution in more generality:
	\begin{mini*}|s|
		{\b f\in \C^0(\X,\R^P)}{\Lcal(\b f)}{}{\bar{\b f} \in}%
		\addConstraint{\b f(x)}{\in \b F(x),\, \forall\, x\in \X.}%
	\end{mini*}%
	For this problem to have a solution, we need to exhibit at least one $\b f\in \C^0(\X,\R^P)$ satisfying $\b f(x)\in \b F(x)$. There are general methods from set-valued analysis to do so, known as \emph{selection theorems} such as the famed \citet[][Theorem 1]{michael1956select}. In \Cref{prop:regular_selection} below, we show that one can even select $\b f\in \C^m(\X,\R^P)$ when allowing to violate the constraints by some small $\epsilon$. The key notion in order to select a function is the \emph{lower semicontinuity} of the set-valued map $\b F$. It is then easy to specialize the notion to affine constraints, for which $\b F(x)$ is a polytope, as in \Cref{lem:lsc_affine constraints} and \Cref{cor:selection_affine_cons}.

	\begin{definition}[Lower semicontinuity of a set-valued map]\label{def:lsc} Let $\X$ be a Hausdorff topological space. A set-valued map $\b F:\X\leadsto \R^P$ is said to be lower semicontinuous at a point $x \in \Dom \b F$ (i.e.\ $\b F(x)\neq \emptyset$) if for any open subset \rv{$\Omega\subset \R^P$} such that $\Omega\cap \b F(x)\neq \emptyset$ there exists a neighborhood $\Ucal$ of $x$ such that 
		\begin{equation}\label{eq:lsc_set-valued}
		\forall\, x'\in \Ucal, \, \Omega\cap \b F(x')\neq \emptyset.
		\end{equation}
		We say that $\b F$ is lower semicontinuous if it is lower semicontinuous at every $x\in\X$.
	\end{definition}
	An equivalent definition for metric spaces $\X$ is that for any sequence $(x_n)_{n\in\N}$ converging to $x$ and $\b y\in \b F(x)$, there is a sequence $\b y_n\in\b  F(x_n)$ converging to $\b y$. It is not always easy to check that a set-valued map is lower semicontinuous. A sufficient criterion for $\b F$ to be lower semicontinuous is that, for any $\epsilon>0$, there exists a neighborhood $\Ucal$ of $x$ such that:
	\begin{equation}\label{eq:eps-lsc_set-valued}
	\forall\, x'\in \Ucal, \, \b F(x)\subset \b F(x')+\epsilon \BB_{\R^P}.
	\end{equation}
	This statement is equivalent to lower semicontinuity if $\b F(x)$ is a compact set \citep[pp43-45]{aubin1984diff}. For instance it allows to check the lower semicontinuity of maps defined by inequalities.
	\begin{lemma}
	Fix $\b g\in C^{0,1}(\X\times \R^P,\R^I)$ such that for any $(x,\b y)$ satisfying $\b g(x,\b y)\le \b 0$, component-wise, there exists a unit-norm $\b v\in\R^P$ for which $Jac_y \b g(x,y)\b v < \b 0$. Then $\b F(x):=\{\b y\,|\,\b g(x,\b y)\le \b 0 \}$ is lower semicontinuous if it is non-empty.
	\end{lemma}
	\begin{proof}
		Let $x$, $\b y$ and $\b v$ be as in the assumptions. Fix $\epsilon_0$ such that: $\forall \epsilon \in (0,\epsilon_0],\, \b g(x,\b y+\epsilon \b v)\b v <  \epsilon \frac{Jac_y \b g(x,y)\b v}{2}$. Let $\epsilon \in (0,\epsilon_0]$ and fix $\eta$ such that $\eta< \epsilon \frac{|Jac_y \b g(x,y)\b v|}{3}$. Let $\Ucal$ be a neighborhood of $x$ such that $\b g(x',\b y+\epsilon \b v)\le \b g(x,\b y+\epsilon \b v) +\eta \b 1$ for $x'\in \Ucal$. Then for all $x'\in \Ucal$, $\b g(x',\b y+\epsilon \b v)\le \epsilon \frac{Jac_y \b g(x,y)\b v}{6}< \b 0 $, hence $\b F$ satisfies \eqref{eq:eps-lsc_set-valued}.
	\end{proof}
	
	\begin{proposition}[Smooth selections]\label{prop:regular_selection}
		Assume that $\X$ is a metric space and $\b F:\X\leadsto \R^P$ is a lower semicontinuous set-valued map with non-empty, closed, convex values. Then there exists $\b f\in \C^0(\X,\R^P)$ such that $\b f(x)\in \b F(x)$ for all $x\in\X$. If, for a given $m\in\N\cup\{\infty\}$, $\X$ is a manifold of class $C^m$, modeled on a separable Hilbert space $\Ecal$, then, for every $\epsilon>0$, there exists $\b f_\epsilon\in \C^m(\X,\R^P)$ such that $\b f_\epsilon(x)\in \b F(x)+\epsilon \BB_{\R^P}(0,1)$ for all $x\in\X$.
	\end{proposition}
	\begin{proof}
		The first statement is only an application of a renowned theorem by Michael \citep[][Theorem 1]{michael1956select} which gives a continuous selection under these assumptions on $F$ for any paracompact space $\X$, metric spaces being paracompact \citep[][Corollary 1]{stone1948paracompact}. Actually $\X$ could just be taken paracompact, while $\R^P$ can be replaced by a Banach \citep{michael1956select} or even Fréchet space \citep{aliprantis2006}.
		
		The second statement is an extension of a result by \citet[][Proof of Theorem 1, p.\ 82]{aubin1984diff} where it was shown for $\b f_\epsilon\in \C^0(\X,\R^P)$ and $\X$ a metric space. We combine the construction of the selection with another result \citep[][Corollary 3.8, p.\ 36]{lang1995diff} on smooth partitions of unity. Since $\b F$ is lower semicontinuous, for every $x\in\X$, we can fix $\b y_x\in \b F(x)$ and $\delta_x>0$ such that $(\b y_x+\epsilon \dBB)\cap \b F(x')\neq \emptyset$ for any $x'\in \BB(x,\delta_x)$. Since $\X$ is a manifold modeled on $\Ecal$, we can also find a chart $(G_x,\gamma_x)$ at every $x\in \X$ (i.e.\ a $C^m$-diffeomorphism $\gamma_x$ from an open neighborhood $G_x$ of $x\in\X$ to an open subset of $\Ecal$) with $G_x\subset \BB(x,\delta_x)$. As any open ball of $\Ecal$ of finite radius is $C^\infty$-isomorphic to $\Ecal$, by composing $\gamma_x$ with such a $C^\infty$-diffeomorphism, without reindexing, we can take $\gamma_x (G_x) = \Ecal$. By definition, $\left\{G_x\right\}_{x\in\X}$ is an open covering of $\X$. Since $\X$ is a metric space, it is paracompact \citep[][Corollary 1]{stone1948paracompact}, so we can find an open refinement\footnote{We can choose for $\left\{U_{x}\right\}_{x\in\X}$ the same index set as $\left\{G_{\X}\right\}_{x\in\X}$  owing to \citep[][Proposition 3.1, p.\ 31]{lang1995diff}.} $\left\{U_{x}\right\}_{x\in\X}$ of the covering $\left\{G_{x}\right\}_{x\in\X}$ (i.e.\ each $U_{x}$ is an open set contained in $G_{x}$ and $\X=\bigcup_{x\in\X}U_{x}$) which is locally finite (i.e.\ for every $x_0\in\X$, there exists a neighborhood of $x_0$ which intersects only a finite number of $U_x$). Let $\varphi_{x}$ be the restriction of $\gamma_{x}$ to $U_{x}$. Again by paracompactness of $\X$, applying \citep[][Proposition 3.3, p.\ 31]{lang1995diff}, we can find two locally finite open refinements of $\left\{U_{x}\right\}_{x\in\X}$, namely $\left\{V_{x}\right\}_{x\in\X}$ and $\left\{W_{x}\right\}_{x\in\X}$, , such that
		$$
		\bar{W}_{x} \subset V_{x} \subset \bar{V}_{x} \subset U_{x}
		$$
		the bar denoting closure in $\X$. Each $\bar{V}_{x}$ being closed in $\X$, as $\varphi_{x}$ is a diffeomorphism, we consequently have that $\varphi_{x} (\bar{V}_{x})$ is closed in $\Ecal$, and so is $\varphi_{x} (\bar{W}_{x})$. Using \citet[][Theorem 3.7, p.\ 36]{lang1995diff}, for each $x\in\X$, we can fix a $C^\infty$-function $\phi_x:\Ecal\rightarrow [0,1]$ such that $\phi_x(\varphi_{x} (\bar{W}_{x}))=\{1\}$ and $\phi_x(\Ecal \backslash\varphi_{x} (V_{x}))=\{0\}$ (which can be constructed based on bump functions). Setting $\psi_{x}=\phi_x \circ \varphi_{x} \in C^m(\X,[0,1])$, $\psi_x$ is equal to $1$ on $\bar{W}_{x}$ and $0$ on $\X\backslash V_{x}$. Thus the mapping $\psi_x:\X\rightarrow \R^P$ defined as 
		\begin{equation*}
		\b f(x'):=\frac{\sum_{x\in\X} \psi_{x}(x') \b y_{x}}{\sum_{x\in\X} \psi_{x}(x')}
		\end{equation*}
		belongs to $C^m(\X,\R^P)$ since it is locally a finite sum of $C^m$ functions (this also justifies the definition of $\b f$ since only a finite number of indexes are active). For any given $x'\in\X$, $\psi_{x}(x')>0$ implies that $x'\in V_x$. However $V_x\subset U_{x}\subset G_x\subset \BB(x,\delta_{x})$. Hence, by definition of $\b y_{x}$, $\b y_{x} \in \b F(x')+\epsilon \dBB$. Since $\b F$ is convex-valued, the convex combination which defines $\b f$ leads to $\b f(x')\in \b F(x')+\epsilon \dBB$, so $\b f(\cdot)$ is the  approximate selection of class $C^m$ which was sought for.
	\end{proof}
	
	\begin{lemma}[l.s.c.\ affine maps]\label{lem:lsc_affine constraints}
		For $\X$ a metric space, $\Kcons$ a non-empty closed subset of $\X$ and $\b  C(\cdot) \in C^0(\Kcons, \R^{I\times P})$, assume that for any $x\in\Kcons$ there exists $\b z_x\in\R^P$ such that $\b C(x)\b z_x>0$ (i.e.\ the constraints are non-degenerate). Then the set-valued map $\b F:x\in\Kcons \leadsto \{\b y\,|\, \b C(x) \b y \ge \b 1\}$ is lower semicontinuous.
	\end{lemma}
	\begin{proof}
		Note that the considered $\b F$ is not in general compact-valued so, rather than rather than \eqref{eq:eps-lsc_set-valued}, we will show that, given $x\in \Kcons$, $\b y\in\R^P$ such that $ \b C(x) \b y > \b 1$, and a sequence $(x_n)_{n\in\N}\in \Kcons^\N$ converging to $x$, there is a sequence $\b y_n\in \b F(x_n)$ converging to $\b y$. Hence the set-valued map $ \mathring{\b F}:x\in\Kcons \leadsto \{\b y\,|\, \b C(x) \b y > \b 1\}$ will be lower semicontinuous by the remark under \Cref{def:lsc}, and we will conclude by a closure argument.
		
		Fix $\b z_x\in\R^P$ such that $\b C(x)\b z_x>0$, thus $\b y_x=3\frac{\b z_x}{\min_i (\b C(x)\b z_x)_i}$ satisfies $\b C(x)\b y_x>3\cdot\b 1$. Since $\b C(x) \b y > \b 1$, we can fix $\epsilon\in(0,1)$ such that $\b C(x) \b y \ge (1+\epsilon)\cdot\b 1$. Then for any $s\in(0,1)$,
		\begin{align*}
		\b C(x_n)(s\b y_x+ (1-s) \b y)&=\b C(x)(s\b y_x+ (1-s) \b y)+(\b C(x_n)- \b C(x))(s\b y_x+ (1-s) \b y)\\
		&\ge (3s+(1-s)(1+\epsilon))\cdot\b 1 +(\b C(x_n)- \b C(x))(s\b y_x+ (1-s) \b y)\\
		&\ge (1+s+\epsilon)\cdot\b 1 -\|\b C(x_n)- \b C(x)\|_2(s\|\b y_x \|_2+(1-s)(\|\b y_x \|_2+\|\b y \|_2)),
		\end{align*}
		whence, by continuity of $\b C(\cdot)$, there exists $N$ independent of $s$ such that \rv{$s\b y_x+ (1-s) \b y  \in \mathring{\b F}(x_n)$} for all $n\ge N$ and all $s\in(0,1) $. Take a sequence $(s_n)_{n\in \N}\in (0,1)^\N$ converging to $0$. Then $\b y_n:=s_n\b y_x+ (1-s_n) \b y\in \mathring{\b F}(x_n)$ and converges to $\b y$, showing that $\mathring{\b F}$ is lower semicontinuous at $x$. To conclude, notice that if in \Cref{def:lsc} we have $\Omega\cap \b F(x)\neq \emptyset$, then $\Omega\cap \mathring{\b F}(x)\neq \emptyset$. Hence \eqref{eq:lsc_set-valued} for $\mathring{\b F}$ gives
  $\Omega\cap \b F(x')\supset \Omega\cap \mathring{\b F}(x') \neq \emptyset$, and $\b F$ is lower semicontinuous at $x$.
		
	\end{proof}
	
	\begin{corollary}[Selection in affine maps]\label{cor:selection_affine_cons}
		Under the assumptions of \Cref{lem:lsc_affine constraints}, if $\b C$ is bounded over $\Kcons$, for every $\zeta\ge 0$, there exists $\b f_\zeta\in C^0 (\X,\R^P)$ such that $\b f_\zeta(x)+\zeta \BB_{\R^P}(0,1) \subset \b F(x)$ for all $x\in\Kcons$. If $\b C$ is unbounded over $\Kcons$, then the result still holds for $\zeta=0$. If, for a given $m\in\N\cup\{\infty\}$, $\X$ is a manifold of class $C^m$, modeled on a separable Hilbert space $\Ecal$, then $\b f_\zeta$ can be taken in $C^m (\X,\R^P)$.
	\end{corollary}
	\begin{proof}
		Fix $\zeta\ge0$, set $M_c:=\max_{x\in\Kcons} \|\b C(x)\|$ and define the set-valued map $\b F_\zeta:x\in\X \leadsto \{\b y\,|\, \b C(x) \b y \ge (1+2\zeta M_c)\b 1\}$. This map satisfies the same assumptions as $\b F$, so $\b F_\zeta$ is also lower semicontinuous over $\Kcons$ and has non-empty, closed, convex values. By \Cref{prop:regular_selection}, there thus exists $\b f_0\in \C^0(\Kcons,\R^P)$ such that $\b f_0(x)\in \b F_\zeta(x)$ for all $x\in\Kcons$. By definition, $\b F_\zeta(x) + 2\zeta \BB_{\R^P}(0,1) \subset \b F(x)$ for all $x\in\Kcons$. So $\b f_0(x)+\zeta \BB_{\R^P}(0,1)\subset \b F(x)$ for all $x\in\Kcons$. Since a metric space is normal and $\Kcons$ is closed, we can apply the Urysohn-Brouwer lemma to continuously extend $\b f_0$ to some $\b f_\zeta$ defined over the whole $\X$. If $M_c=+\infty$, then, for $\zeta=0$, we can still apply \Cref{prop:regular_selection} which yields directly the result.
		
		If $\X$ is more regular, we again apply \Cref{prop:regular_selection} but to the map $\tilde{\b F}(x):=\left\{
		\begin{array}{ll}
		\R^P  \text{ if $x\notin \Kcons$}\\ 
		\b f_0(x)  \text{ if $x\in \Kcons$}
		\end{array} 	\right.	  $. Since $\b f_0\in \C^0(\Kcons,\R^P)$ and $\Kcons$ is closed, the latter is also lower semicontinuous with non-empty, closed, convex values.

	\end{proof}
	
	\bibliographystyle{apalike}
	\bibliography{inequalities,biblioPCAF}
\end{document}